\documentclass{amsart}

\usepackage{graphicx}        
\usepackage{multicol}        
\usepackage[bottom]{footmisc}
\usepackage{times}

\usepackage{amsmath}
\usepackage{amssymb}

\usepackage{hyperref}

\newtheorem{theorem}{Theorem}[section]

\newtheorem{proposition}[theorem]{Proposition}
\newtheorem{corollary}[theorem]{Corollary}
\theoremstyle{definition}
\newtheorem{definition}[theorem]{Definition}
\newtheorem{example}[theorem]{Example}
\newtheorem{remark}[theorem]{Remark}

\makeindex             

\usepackage{maplestd2e}

\begin{document}

\title{Thomas Decomposition and Nonlinear Control Systems}
\author{Markus Lange-Hegermann and Daniel Robertz}
%
%
\maketitle

\begin{abstract}
This paper applies the Thomas decomposition technique to nonlinear control
systems, in particular to the study of the dependence of the system behavior on
parameters. Thomas' algorithm is a symbolic method which splits
a given system of nonlinear partial differential equations into a finite family
of so-called simple systems which are formally integrable and define a partition
of the solution set of the original differential system. Different simple systems
of a Thomas decomposition describe different structural behavior of the control
system in general. The paper gives an introduction to the Thomas decomposition
method and shows how notions such as invertibility, observability and flat outputs
can be studied. A Maple implementation of Thomas' algorithm is used to illustrate
the techniques on explicit examples.
\end{abstract}

\section{Introduction}
\label{mlhdr:sec:introduction}

This paper gives an introduction to the Thomas decomposition method and
presents first steps in applying it to the structural study of nonlinear
control systems. It extends and refines our earlier work \cite{LangeHegermannRobertz1}.

\medskip

Symbolic computation allows to study many structural aspects
of control systems, e.g., controllability, observability,
input-output behavior, etc. In contrast to a numerical
treatment, the dependence of the results on parameters occurring
in the system is accessible to symbolic methods.

\medskip

An algebraic approach for treating nonlinear control systems
has been developed during the last decades, e.g., by M.~Fliess
and coworkers, J.-F.\ Pommaret and others, cf., e.g., \cite{FliessGlad}, \cite{Glad},
\cite{Pommaret2001}, and the references therein. In particular, the notion
of flatness has been studied extensively and has been applied to many
interesting control problems (cf.,
e.g., \cite{FliessLevineMartinRouchon}, \cite{DavidJeanBaptiste}, \cite{Levine11}).
The approach of Diop \cite{DiopDiffAlg,DiopElimination} builds on the
characteristic set method (cf.\ \cite{Kolchin}, \cite{WuMathMech}).
The Rosenfeld-Gr{\"o}bner algorithm (cf.\ \cite{BoulierLazardOllivierPetitotAAECC})
can be used to perform the relevant computations effectively; implementations
of related techniques are available, e.g., as Maple packages
{\tt DifferentialAlgebra} (by F.~Boulier and E.~S.\ Cheb-Terrab),
formerly {\tt diffalg} (by F.~Boulier and E.~Hubert), and
{\tt RegularChains} (by F.~Lemaire, M.~Moreno Maza, and Y.\ Xie)
\cite{RegularChains}; cf.\ also \cite{WangEliminationMethods}
for alternative approaches.
As an example of an application of the Rosenfeld-Gr{\"o}bner algorithm
we refer to \cite{PicoMarco}, where it is demonstrated how to compute
a block feedforward form and a generalized controller form for a
nonlinear control system.

\medskip

So far the dependence of nonlinear control systems on parameters has not
been studied by a rigorous method such as Thomas decomposition.
This paper demonstrates how the Thomas decomposition method can be
applied in this context. In particular, Thomas' algorithm can detect certain
structural properties of control systems by performing elimination
and it can separate singular cases of behavior in control systems from
the generic case due to splitting into disjoint solution sets. We also
consider the Thomas decomposition method as a preprocessing technique
for the study of a linearization of a nonlinear system
(cf.\ \cite[Sect.~5.5]{Robertz7}), an aspect that we do not pursue here.

\medskip

Dependence of control systems on parameters has been examined, in
particular, by J.-F.\ Pommaret and A.~Quadrat
in \cite{PommaretQuadrat1997}, \cite{Pommaret2001}.
For linear systems, stratifications of the space of
parameter values have been studied using Gr{\"o}bner bases
in \cite{LevandovskyyZerz}.

\medskip

In the 1930s the American mathematician J.~M.\ Thomas designed an
algorithm which decomposes a polynomially nonlinear system of
partial differential equations into so-called simple systems.
The algorithm uses, in contrast to the characteristic set method,
inequations to provide a disjoint decomposition of the solution
set (cf.\ \cite{Thomas}).
It precedes work by E.~R.\ Kolchin \cite{Kolchin} and A.~Seidenberg \cite{Seidenberg},
who followed J.~F.\ Ritt \cite{Ritt}.
Recently a new algorithmic approach to the Thomas decomposition method
has been developed (cf.\
\cite{GerdtThomas,BaechlerGerdtLangeHegermannRobertz2,Robertz6}),
building also on ideas of the French mathematicians C.~Riquier \cite{Riquier}
and M.~Janet \cite{Janet1}.
Implementations as Maple packages of the algebraic and differential parts
of Thomas' algorithm are available due to work by T.~B{\"a}chler and M.~Lange-Hegermann
\cite{ThomasPackages}.
The implementation of the differential part is available in the
Computer Physics Communications library \cite{GLHR19}
and has also been incorporated into Maple's standard library since Maple 2018.
An earlier implementation of the algebraic part was given by D.~Wang
\cite{WangSimpleSystems}.

\medskip

Section~\ref{mlhdr:sec:thomas} introduces the Thomas decomposition method
for algebraic and differential systems and discusses the main properties
of its output. The algorithm for the differential case builds on the
algebraic part. Section~\ref{mlhdr:sec:elimination} explains how the
Thomas decomposition technique can be used to solve elimination problems
that occur in our study of nonlinear control systems. Finally,
Section~\ref{mlhdr:sec:control} addresses concepts of nonlinear control
theory, such as invertibility, observability, and flat outputs, possibly
depending on parameters of the control system, and gives examples using
a Maple implementation of Thomas' algorithm.

\section{Thomas decomposition}
\label{mlhdr:sec:thomas}

This section gives an introduction to the Thomas decomposition method
for algebraic and differential systems. The case of differential systems,
discussed in Subsection~\ref{mlhdr:subsec:diffsystem}, builds on the case
of algebraic systems which is dealt with in the first subsection.
For more details on Thomas' algorithm, we refer
to \cite{BaechlerGerdtLangeHegermannRobertz2}, \cite{GerdtThomas},
\cite{PleskenCounting}, \cite{BaechlerThesis}, \cite{LangeHegermannThesis},
and \cite[Sect.~2.2]{Robertz6}.

\subsection{Algebraic systems}
\label{mlhdr:subsec:algsystem}

Let $K$ be a field of characteristic zero and $R = K[x_1, \ldots, x_n]$
the polynomial algebra with indeterminates $x_1$, \ldots, $x_n$ over $K$.
We denote by $\overline{K}$ an algebraic closure of $K$.

\begin{definition}
An \emph{algebraic system} $S$, defined over $R$, is given by finitely
many equations and inequations
\begin{equation}\label{mlhdr:eq:algsys}
p_1 = 0, \quad p_2 = 0, \quad \ldots, \quad p_s = 0, \quad
q_1 \neq 0, \quad q_2 \neq 0, \quad \ldots, \quad q_t \neq 0,
\end{equation}
where $p_1$, \ldots, $p_s$, $q_1$, \ldots, $q_t \in R$ and
$s$, $t \in \mathbb{Z}_{\ge 0}$.
The \emph{solution set} of $S$ in $\overline{K}^n$ is
\[
{\rm Sol}_{\overline{K}}(S) := \{ \, a \in \overline{K}^n \mid
p_i(a) = 0 \mbox{ and } q_j(a) \neq 0 \mbox{ for all }
1 \le i \le s, \, 1 \le j \le t \, \}.
\]
\end{definition}

\smallskip

We fix a total ordering $>$ on the set $\{ x_1, \ldots, x_n \}$ allowing
us to consider every non-constant element $p$ of $R$ as a univariate
polynomial in the greatest variable with respect to $>$ which occurs in $p$,
with coefficients which are themselves univariate polynomials in lower
ranked variables, etc. Without loss of generality we may assume that
$x_1 > x_2 > \ldots > x_n$. The choice of $>$ corresponds to a choice of
projections
\[
\begin{array}{rclrcr}
\pi_1\colon \overline{K}^n & \longrightarrow & \overline{K}^{n-1}\colon &
(a_1, a_2, \ldots, a_n) & \longmapsto & (a_2, a_3, a_4, \ldots, a_n),\\[0.5em]
\pi_2\colon \overline{K}^n & \longrightarrow & \overline{K}^{n-2}\colon &
(a_1, a_2, \ldots, a_n) & \longmapsto & (a_3, a_4, \ldots, a_n),\\[0.5em]
& \vdots & & & \vdots \\[0.5em]
\pi_{n-1}\colon \overline{K}^n & \longrightarrow & \overline{K}\colon &
(a_1, a_2, \ldots, a_n) & \longmapsto & a_n.
\end{array}
\]
Thus, the recursive representation of polynomials is motivated by
considering each $\pi_{k-1}({\rm Sol}_{\overline{K}}(S))$ as fibered
over $\pi_k({\rm Sol}_{\overline{K}}(S))$, for $k = 1$, \ldots, $n-1$,
where $\pi_0 := {\rm id}_{\overline{K}^n}$ (cf.\ also \cite{PleskenCounting}).
The purpose of a Thomas decomposition of ${\rm Sol}_{\overline{K}}(S)$,
to be defined below, is to clarify this fibration structure. The
solution set ${\rm Sol}_{\overline{K}}(S)$ is partitioned into subsets
${\rm Sol}_{\overline{K}}(S_1)$, \ldots, ${\rm Sol}_{\overline{K}}(S_r)$
in such a way that, for each $i = 1$, \ldots, $r$ and $k = 1$, \ldots, $n-1$,
the fiber cardinality $|\pi_k^{-1}(\{ \, a \, \})|$
does not depend on the choice of $a \in \pi_k({\rm Sol}_{\overline{K}}(S_i))$.
In terms of the defining equations and inequations in (\ref{mlhdr:eq:algsys}),
the fundamental obstructions to this uniform behavior are zeros of the
leading coefficients of $p_i$ or $q_j$ and zeros of $p_i$ or $q_j$ of
multiplicity greater than one.

\begin{definition}\label{mlhdr:de:prep}
\renewcommand{\theenumi}{\alph{enumi})}
\renewcommand{\labelenumi}{\theenumi}
\renewcommand{\theenumii}{\roman{enumii}}
Let $p \in R \setminus K$.
\begin{enumerate}
\item \label{mlhdr:leader}
The greatest variable with respect to $>$ which occurs in $p$ is referred
to as the \emph{leader} of $p$ and is denoted by ${\rm ld}(p)$.
\item \label{mlhdr:degree}
For $v = {\rm ld}(p)$ we denote by $\deg_{v}(p)$ the degree of $p$ in $v$.
\item \label{mlhdr:initial}
The coefficient of the highest power of ${\rm ld}(p)$ occurring in $p$ is
called the \emph{initial} of $p$ and is denoted by ${\rm init}(p)$.
\item \label{mlhdr:discriminant}
The \emph{discriminant} of $p$ is defined as
\[
{\rm disc}(p) := (-1)^{d(d-1)/2} \, {\rm res}\left( p,
\frac{\partial p}{\partial \, {\rm ld}(p)}, {\rm ld}(p) \right) / \,
{\rm init}(p), \qquad
d = \deg_{{\rm ld}(p)}(p),
\]
where ${\rm res}(p, q, v)$ is the resultant of $p$ and $q$ with respect to
the variable $v$.
(Note that ${\rm disc}(p)$ is a polynomial because ${\rm init}(p)$ divides
${\rm res}(p, \partial p/\partial \, {\rm ld}(p), {\rm ld}(p))$,
since the Sylvester matrix, whose
determinant is ${\rm res}(p, \partial p/\partial \, {\rm ld}(p), {\rm ld}(p))$,
has a column all of whose entries are divisible by ${\rm init}(p)$.)
\end{enumerate}
\end{definition}

Both ${\rm init}(p)$ and ${\rm disc}(p)$ are elements of the
polynomial algebra $K[x \mid x < {\rm ld}(p)]$. The zeros of a
univariate polynomial which have multiplicity greater than one are
the common zeros of the polynomial and its derivative. The solutions
of ${\rm disc}(p) = 0$ in $\overline{K}^{n-k}$, where ${\rm ld}(p) = x_k$,
are therefore those tuples $(a_{k+1}, a_{k+2}, \ldots, a_n)$ for which the
substitution $x_{k+1} = a_{k+1}$, $x_{k+2} = a_{k+2}$, \ldots, $x_n = a_n$
in $p$ results in a univariate polynomial with a zero of multiplicity greater than one.

\begin{definition}\label{mlhdr:de:algsimple}
\renewcommand{\theenumi}{\alph{enumi})}
\renewcommand{\labelenumi}{\theenumi}
\renewcommand{\theenumii}{\roman{enumii}}
An algebraic system $S$, defined over $R$, as in (\ref{mlhdr:eq:algsys})
is said to be \emph{simple}
(with respect to $>$) if the following three conditions hold.
\begin{enumerate}
\item \label{mlhdr:algsimple1}
For all $i = 1$, \ldots, $s$ and $j = 1$, \ldots, $t$ we have $p_i \not\in K$ and $q_j \not\in K$.
\item \label{mlhdr:algsimple2}
The leaders of the left hand sides of the equations and inequations
in $S$ are pairwise different, i.e.,
$|\{ \, {\rm ld}(p_1), \ldots, {\rm ld}(p_s), {\rm ld}(q_1), \ldots, {\rm ld}(q_t) \, \}| = s+t$.
\item \label{mlhdr:algsimple3}
For every $r \in \{ \, p_1, \ldots, p_s, q_1, \ldots, q_t \, \}$, if
${\rm ld}(r) = x_k$, then neither of the equations
${\rm init}(r) = 0$ and ${\rm disc}(r) = 0$
has a solution $(a_{k+1}, a_{k+2}, \ldots, a_n)$ in $\pi_k({\rm Sol}_{\overline{K}}(S))$.
\end{enumerate}
\end{definition}

Subsets of non-constant polynomials in $R$ with pairwise different leaders
(i.e., satisfying \ref{mlhdr:algsimple1} and \ref{mlhdr:algsimple2})
are also referred to as triangular sets (cf., e.g.,
\cite{AubryLazardMorenoMaza}, \cite{Hubert},
\cite{WangEliminationMethods}).

\begin{remark}
A simple algebraic system $S$ admits the following solution procedure, which
also shows that its solution set is not empty.
Let $S_{< k}$ be the subset of $S$ consisting of the
equations $p = 0$ and inequations $q \neq 0$ with ${\rm ld}(p) < x_k$
and ${\rm ld}(q) < x_k$.
The fibration structure implied by \ref{mlhdr:algsimple3} ensures that,
for every $k = 1$, \ldots, $n-1$, every solution $(a_{k+1}, a_{k+2}, \ldots, a_n)$
of $\pi_k({\rm Sol}_{\overline{K}}(S)) = \pi_k({\rm Sol}_{\overline{K}}(S_{< k}))$
can be extended to a solution $(a_k, a_{k+1}, \ldots, a_n)$
of $\pi_{k-1}({\rm Sol}_{\overline{K}}(S))$.
If $S$ contains an equation $p = 0$ with leader $x_k$,
then there exist exactly $\deg_{x_k}(p)$ such elements $a_k \in \overline{K}$
(because zeros with multiplicity greater than one are excluded by the
non-vanishing discriminant).
If $S$ contains an inequation $q \neq 0$ with leader $x_k$,
all $a_k \in \overline{K}$ except $\deg_{x_k}(q)$ elements define a
tuple $(a_k, a_{k+1}, \ldots, a_n)$ as above. If no equation and no
inequation in $S$ has leader $x_k$, then $a_k \in \overline{K}$ can be
chosen arbitrarily.
\end{remark}

\begin{definition}
Let $S$ be an algebraic system, defined over $R$.
A \emph{Thomas decomposition} of $S$ (or of ${\rm Sol}_{\overline{K}}(S)$)
with respect to $>$ is a collection of finitely many simple
algebraic systems $S_1$, \ldots, $S_r$, defined over $R$, such that
${\rm Sol}_{\overline{K}}(S)$ is the disjoint union of the solution
sets ${\rm Sol}_{\overline{K}}(S_1)$, \ldots, ${\rm Sol}_{\overline{K}}(S_r)$.
\end{definition}

We outline Thomas' algorithm for computing a Thomas decomposition of
algebraic systems.

\begin{remark}\label{mlhdr:rem:algebraicpart}
Given $S$ as in (\ref{mlhdr:eq:algsys}) and a total ordering $>$
on $\{ x_1, \ldots, x_n \}$, a Thomas decomposition of $S$ with respect to $>$
can be constructed by combining Euclid's algorithm with a splitting strategy.

First of all, if $S$ contains an equation $c = 0$ with $0 \neq c \in K$
or the inequation $0 \neq 0$, then $S$ is discarded because it has no solutions.
Moreover, from now on the equation $0 = 0$ and inequations $c \neq 0$
with $0 \neq c \in K$ are supposed to be removed from $S$.

An elementary step of the algorithm
applies a pseu\-do-di\-vi\-sion to a pair $p_1$, $p_2$ of
non-constant polynomials in $R$ with the same leader $x_k$ and
$\deg_{x_k}(p_1) \ge \deg_{x_k}(p_2)$.
The result is a pseu\-do-re\-main\-der
\begin{equation}\label{mlhdr:eq:pseudoreduction}
r = c_1 \cdot p_1 - c_2 \cdot p_2,
\end{equation}
where $c_1$, $c_2 \in R$ and $r$ is constant or has leader less than $x_k$ or
has leader $x_k$ and $\deg_{x_k}(r) < \deg_{x_k}(p_1)$. Since the coefficients
of $p_1$ and $p_2$ are polynomials in lower ranked variables, multiplication
of $p_1$ by a non-constant polynomial $c_1$ may be necessary in general to perform
the reduction in $R$ (and not in its field of fractions). The choice of $c_1$ as a
suitable power of ${\rm init}(p_2)$ always achieves this.

In order to turn $S$ into a triangular set, the algorithm deals with three
kinds of subsets of $S$ of cardinality two. Firstly, each pair of equations
$p_1 = 0$, $p_2 = 0$ in $S$ with ${\rm ld}(p_1) = {\rm ld}(p_2)$ is replaced
with the single equation $r = 0$, where $r$ is the result of applying
Euclid's algorithm to $p_1$ and $p_2$, considered as univariate polynomials
in their leader, using the above pseu\-do-di\-vi\-sion. (If this computation was
stable under substitution of values for lower ranked variables in $p_1$ and $p_2$,
then $r$ would be the greatest common divisor of the specialized polynomials.)

The solution set of the system is supposed not to change, when the
equation $p_1 = 0$ is replaced with the equation $r = 0$ given by the
pseu\-do-re\-duction (\ref{mlhdr:eq:pseudoreduction}). Therefore, we
assume that the polynomial $c_1$, and hence ${\rm init}(p_2)$, does
not vanish on the solution set of the system. In order to ensure this
condition, a preparatory step splits the system into two, if necessary,
and adds the inequation ${\rm init}(p_2) \neq 0$ to one of them and the
equation ${\rm init}(p_2) = 0$ to the other. The algorithm then deals
with both systems separately. These case distinctions also allow to
arrange for the part of condition~\ref{mlhdr:algsimple3} in Definition~\ref{mlhdr:de:algsimple} which concerns initials.

Secondly, let $p = 0$, $q \neq 0$ be in $S$
with ${\rm ld}(p) = {\rm ld}(q) = x_k$. If $\deg_{x_k}(p) \le \deg_{x_k}(q)$,
then $q \neq 0$ is replaced with $r \neq 0$, where $r$ is the result of
applying the pseu\-do-di\-vi\-sion (\ref{mlhdr:eq:pseudoreduction}) to $q$ and $p$.
Otherwise, Euclid's algorithm is applied to $p$ and $q$, keeping track of
the coefficients used for the reductions as in (\ref{mlhdr:eq:pseudoreduction}).
Given the result $r$,
the system is then split into two, adding the conditions $r \neq 0$ and $r = 0$,
respectively. The
inequation $q \neq 0$ is removed from the first new system, because $p = 0$
and $q \neq 0$ have no common solution in that case.
The assumption $r = 0$ and the bookkeeping allows to divide $p$ by the
common factor of $p$ and $q$ (modulo left hand sides of equations with
smaller leader). The left hand side of $p = 0$ is replaced with that
quotient in the second new system. Not all of these cases need a closer
inspection. For instance, if $p$ divides $q$, then the solution set of $S$
is empty and $S$ is discarded.

Thirdly, for a pair $q_1 \neq 0$, $q_2 \neq 0$ in $S$
with ${\rm ld}(q_1) = {\rm ld}(q_2)$, Euclid's algorithm is applied to $q_1$
and $q_2$ in the same way as above. Keeping track of the coefficients
used in intermediate steps allows to determine the least common multiple $m$
of $q_1$ and $q_2$, which again depends on distinguishing the cases whether
the result of Euclid's algorithm vanishes or not. The pair $q_1 \neq 0$, $q_2 \neq 0$
is then replaced with $m \neq 0$.

The part of condition~\ref{mlhdr:algsimple3} in Definition~\ref{mlhdr:de:algsimple}
regarding discriminants is taken care of by applying Euclid's algorithm as
above to $p$ and $\partial p/\partial \, {\rm ld}(p)$, where $p$ is the
left hand side of an equation or inequation. Bookkeeping allows to determine
the square-free part of $p$, which depends again on case distinctions.

Expressions tend to grow very quickly when performing these reductions, so
that an appropriate strategy is essential for dealing with non-trivial systems.
Apart from dividing by the content (in $K$) of polynomials, in intermediate steps
of Euclid's algorithm the coefficients should be reduced modulo equations
in the system with lower ranked leaders. In practice, subresultant computations
(cf., e.g., \cite{Mishra}) allow to diminish the growth of coefficients significantly.

Termination of the procedure sketched above depends on the organization of
its steps. One possible strategy is to maintain an intermediate triangular
set, reduce new equations and inequations modulo the equations in the
triangular set, and select among these results the one with smallest leader
and least degree, preferably an equation, for insertion into the triangular
set. If the set already contains an equation or inequation with the same
leader, then the pair is treated as discussed above. Since equations are
replaced with equations of smaller degree and inequations are replaced with
equations if possible or with the least common multiple of inequations, this
strategy terminates after finitely many steps.

For more details on the algebraic part of Thomas' algorithm,
we refer to \cite{BaechlerGerdtLangeHegermannRobertz2}, \cite{BaechlerThesis},
and \cite[Subsect.~2.2.1]{Robertz6}.

An implementation of Thomas' algorithm for algebraic systems has been
developed by T.~B{\"a}chler as Maple package {\tt AlgebraicThomas}
\cite{ThomasPackages}.
\end{remark}

In what follows, variables are underlined to emphasize that they
are leaders of polynomials with respect to the fixed total ordering $>$.

\begin{example}\label{mlhdr:ex:thomascircle}
Let us compute a Thomas decomposition of the algebraic system
\[
x^2 + y^2 - 1 = 0
\]
consisting of one equation, defined over $R = \mathbb{Q}[x, y]$,
with respect to $x > y$. We set $p_1 := x^2+y^2-1$. Then we have
${\rm ld}(p_1) = x$ and ${\rm init}(p_1) = 1$ and
\[
{\rm disc}(p_1) = -4 \, y^2 + 4.
\]
We distinguish the cases whether or not $p_1 = 0$ has a solution
which is also a zero of ${\rm disc}(p_1)$, or equivalently, of
$y^2 - 1$. In other words, we replace the original algebraic system
with two algebraic systems which are obtained by adding the
inequation $y^2 - 1 \neq 0$ or the equation $y^2 - 1 = 0$. The
first system is readily seen to be simple, whereas the second one
is transformed into a simple system by taking the difference of
the two equations and computing a square-free part. Clearly, the
solution sets of the two resulting simple systems form a partition
of the solution set of $p_1 = 0$. We obtain the Thomas decomposition
\begin{center}
\fbox{\;
$\begin{array}{rcl}
\underline{x}^2 + y^2 - 1 & = & 0\\[1em]
\underline{y}^2 - 1 & \neq & 0
\end{array}$
\;}
\qquad
\fbox{\;
$\begin{array}{rcl}
\underline{x} & = & 0\\[1em]
\underline{y}^2 - 1 & = & 0
\end{array}$
\;}
\end{center}
In this example, all points of ${\rm Sol}_{\overline{K}}(\{ \, p_1 = 0 \, \})$
for which the projection $\pi_1$ onto the $y$-axis has fibers of an
exceptional cardinality have real coordinates, and the significance
of the above case distinction can be confirmed graphically.

As a further illustration let us augment the original system by the
equation which expresses the coordinate $t$ of the point of intersection
of the line through the two points $(0, 1)$ and $(x, y)$ on the circle
with the $x$-axis (stereographic projection): 
\[
\left\{
\begin{array}{rcl}
x^2 + y^2 - 1 & = & 0\\[1em]
(1-y) \, t - x & = & 0
\end{array}
\right.
\]
A Thomas decomposition with respect to $x > y > t$ is obtained as
follows. We set $p_2 := x + t \, y - t$. Since ${\rm ld}(p_1) = {\rm ld}(p_2)$,
we apply polynomial division:
\[
p_1 - (x-t \, y+t) \, p_2 = (1+t^2) \, \underline{y}^2 -
2 \, t^2 \, \underline{y} + t^2 - 1 =
(\underline{y} - 1) \, ((1+t^2) \, \underline{y} - t^2 + 1)\,.
\]
Replacing $p_1$ with the remainder of this division does not alter
the solution set of the algebraic system. It is convenient (but not
necessary) to split the system into two systems according to the
factorization of the remainder:
\[
\left\{ \begin{array}{rcl}
\underline{x} + t \, y - t & = & 0\\[1em]
(1+t^2) \, \underline{y} - t^2 + 1 & = & 0\\[1em]
\underline{y} - 1 & \neq & 0
\end{array} \right.
\qquad \qquad
\left\{ \begin{array}{rcl}
\underline{x} + t \, y - t & = & 0\\[1em]
\phantom{x}\\[1em]
\underline{y} - 1 & = & 0
\end{array} \right.
\]
Another polynomial division reveals that the equation and the inequation
with leader $y$ in the first system have no common solutions. Therefore,
the inequation can be omitted from that system. The initial of the equation
has to be investigated. In fact, the assumption $1+t^2=0$ leads to a
contradiction. Finally, the equation with leader $y$ can be used to
eliminate $y$ in the equation with leader $x$:
\[
(1+t^2) \, (\underline{x} + t \, y - t) - t \,
((1+t^2) \, \underline{y} - t^2 + 1) = (1+t^2) \, \underline{x} - 2 \, t\,.
\]
A similar simplification can be applied to the second system above. We
obtain the Thomas decomposition
\begin{center}
\fbox{\;
$\begin{array}{rcl}
(1+t^2) \, \underline{x} - 2 \, t & = & 0\\[1em]
(1+t^2) \, \underline{y} - t^2 + 1 & = & 0\\[1em]
\underline{t}^2 + 1 & \neq & 0
\end{array}$
\;}
\qquad
\fbox{\;
$\begin{array}{rcl}
\underline{x} & = & 0\\[1em]
\underline{y} - 1 & = & 0\\[1em]
\phantom{x}
\end{array}$
\;}
\end{center}
from which a rational parametrization of the circle can be read off.
\end{example}

\begin{remark}\label{mlhdr:rem:thomasdecompnonunique}
A Thomas decomposition of an algebraic system is not uniquely determined.
It depends on the chosen total ordering $>$, the order in which
intermediate systems are dealt with and other choices, such as whether
factorizations of left hand sides of equations are taken into account or not.
\end{remark}

According to Hilbert's Nullstellensatz (cf., e.g., \cite{Eisenbud}),
the solution sets $V$ in $\overline{K}^n$ of systems of polynomial equations in $x_1$, \ldots, $x_n$,
defined over $R$, are in one-to-one correspondence with their
vanishing ideals in $R$
\[
\mathcal{I}_R(V) := \{ \, p \in R \mid p(a) = 0 \mbox{ for all } a \in V \, \},
\]
and these are the radical ideals of $R$, i.e., the
ideals $I$ of $R$ which equal their radicals
\[
\sqrt{I} := \{ \, p \in R \mid p^r \in I \mbox{ for some } r \in \mathbb{Z}_{\ge 0} \, \}.
\]
The solution sets $V$ can then be considered as the closed subsets
of $\overline{K}^n$ with respect to the Zariski topology.

The fibration structure of a simple algebraic system $S$ allows to deduce
that the polynomials in $R$ which vanish on ${\rm Sol}_{\overline{K}}(S)$
are precisely those polynomials in $R$ whose pseu\-do-re\-main\-ders modulo $p_1$,
\ldots, $p_s$ are zero, where $p_1 = 0$, \ldots, $p_s = 0$ are the
equations in $S$. If $E$ is the ideal of $R$ generated by $p_1$, \ldots, $p_s$
and $q$ the product of all ${\rm init}(p_i)$, then these polynomials form
the saturation ideal
\[
E : q^{\infty} := \{ \, p \in R \mid q^r \cdot p \in E \mbox{ for some } r \in \mathbb{Z}_{\ge 0} \, \}.
\]
In particular, simple algebraic systems admit an effective way to
decide membership of a polynomial to the associated radical ideal
(cf.\ also Proposition~\ref{mlhdr:prop:intersection} below).

\begin{proposition}[\cite{Robertz6}, Prop.~2.2.7]\label{mlhdr:prop:algebraicmembership}
Let $S$ be a simple algebraic system as in (\ref{mlhdr:eq:algsys}),
$E$ the ideal of $R$ generated by $p_1$, \ldots, $p_s$,
and $q$ the product of all ${\rm init}(p_i)$.
Then $E : q^{\infty}$ consists of all polynomials
in $R$ which vanish on ${\rm Sol}_{\overline{K}}(S)$. In particular,
$E : q^{\infty}$ is a radical ideal.
Given $p \in R$, we have $p \in E : q^{\infty}$ if and only if the
pseu\-do-re\-main\-der of $p$ modulo $p_1$, \ldots, $p_s$ is zero.
\end{proposition}

\subsection{Differential systems}
\label{mlhdr:subsec:diffsystem}

\begin{definition}
A \emph{differential field} $K$
with commuting derivations $\delta_1$, \ldots, $\delta_n$
is a field $K$ endowed with maps $\delta_i\colon K \to K$,
satisfying
\[
\delta_i(k_1 + k_2) = \delta_i(k_1) + \delta_i(k_2), \quad
\delta_i(k_1 \, k_2) = \delta_i(k_1) \, k_2 + k_1 \, \delta_i(k_2) \quad
\mbox{for all } k_1, k_2 \in K,
\]
$i = 1$, \ldots, $n$, and
$\delta_i \circ \delta_j = \delta_j \circ \delta_i$ for all $1 \le i, j \le n$.
\end{definition}

In what follows, let $K$ be the differential field of (complex)
meromorphic functions on an open and connected subset $\Omega$
of $\mathbb{C}^n$. The derivations on $K$ are given by the partial
differential operators $\delta_1$, \ldots, $\delta_n$ with respect
to the coordinates of $\mathbb{C}^n$. Moreover, let $R = K\{ u_1, \ldots, u_m \}$
be the differential polynomial ring in the
differential indeterminates $u_1$, \ldots, $u_m$.
These indeterminates give rise to symbols $(u_k)_J$,
where $J = (j_1, \ldots, j_n) \in (\mathbb{Z}_{\ge 0})^n$, which represent
the partial derivatives of $m$ infinitely differentiable functions.
More precisely, $R$ is the polynomial algebra $K[(u_k)_J \mid 1 \le k \le m, \, J \in (\mathbb{Z}_{\ge 0})^n]$ over $K$ in infinitely many indeterminates $(u_k)_J$,
endowed with commuting derivations $\partial_1$, \ldots, $\partial_n$ such that
\[
\partial_j \, ((u_k)_J) = (u_k)_{J+1_j}, \qquad
\partial_j|_K = \delta_j \quad \mbox{for all } j = 1, \ldots, n,
\]
where $1_j$ is the $j$-th standard basis vector
of $\mathbb{Z}^n$. For $k \in \{ 1, \ldots, m \}$,
we identify $(u_k)_{(0, \ldots, 0)}$ and $u_k$.
We set $\Delta := \{ \, \partial_1, \ldots, \partial_n \, \}$, and for any
subset $\{ \, \partial_{i_1}, \ldots, \partial_{i_r} \, \}$ of $\Delta$
we define the free commutative monoid
of all monomials in $\partial_{i_1}$, \ldots, $\partial_{i_r}$
\[
{\rm Mon}(\{ \, \partial_{i_1}, \ldots, \partial_{i_r} \, \}) :=
\{ \, \partial_{i_1}^{e_1} \ldots \partial_{i_r}^{e_r} \mid e \in (\mathbb{Z}_{\ge 0})^r \, \}.
\]

\begin{definition}
A \emph{differential system} $S$, defined over $R = K\{ u_1, \ldots, u_m \}$, is
given by finitely many equations and inequations
\begin{equation}\label{mlhdr:eq:diffsys}
p_1 = 0, \quad p_2 = 0, \quad \ldots, \quad p_s = 0, \quad
q_1 \neq 0, \quad q_2 \neq 0, \quad \ldots, \quad q_t \neq 0,
\end{equation}
where $p_1$, \ldots, $p_s$, $q_1$, \ldots, $q_t \in R$ and
$s$, $t \in \mathbb{Z}_{\ge 0}$.
The \emph{solution set} of $S$ is
\begin{eqnarray*}
{\rm Sol}_{\Omega}(S) := \{ \, f = (f_1, \ldots, f_m) & | & f_k\colon \Omega \to \mathbb{C} \mbox{ analytic}, \, k = 1, \ldots, m, \\[0.2em]
& & p_i(f) = 0, \, q_j(f) \neq 0, \, i = 1, \ldots, s, \, j = 1, \ldots, t \, \}.
\end{eqnarray*}
\end{definition}

\begin{remark}
Since each component $f_k$ of a solution of (\ref{mlhdr:eq:diffsys}) is
assumed to be analytic, the equations $p_i = 0$ and inequations $q_j \neq 0$ 
(and their consequences) can be translated into algebraic conditions on the
Taylor coefficients of power series expansions of $f_1$, \ldots, $f_m$
(around a point in $\Omega$). An inequation $q \neq 0$ then turns into a
disjunction of algebraic inequations for all coefficients which result
from substitution of power series expansions for $u_1$, \ldots, $u_m$ in $q$.
(This approach leads to the definition of the differential counting polynomial,
a fine invariant of a differential system \cite{LH18}).

An appropriate choice of $\Omega \subseteq \mathbb{C}^n$ can often only be
made after the formal treatment of a given differential system by Thomas'
algorithm (as, e.g., singularities of coefficients in differential consequences
will only be detected during that process). In general, we assume that $\Omega$
is chosen in such a way that the given systems have analytic solutions on $\Omega$.
\end{remark}

Clearly, by neglecting the derivations on $R = K\{ u_1, \ldots, u_m \}$, a
differential system can be considered as an algebraic system in the
finitely many variables $(u_i)_J$ which occur in the equations and inequations.
The same recursive representation of polynomials as in the algebraic
case is employed, but the total ordering on the set of variables $(u_i)_J$
is supposed to respect the action of the derivations.

\begin{definition}
\renewcommand{\theenumi}{\alph{enumi})}
\renewcommand{\labelenumi}{\theenumi}
\renewcommand{\theenumii}{\roman{enumii}}
A \emph{ranking} $>$ on $R = K\{ u_1, \ldots, u_m \}$ is a total ordering on the set
\[
{\rm Mon}(\Delta) \, u := \{ \, (u_k)_J \mid 1 \le k \le m, \, J \in (\mathbb{Z}_{\ge 0})^n \, \}
\]
such that for all $j \in \{ 1, \ldots, n \}$,
$k$, $k_1$, $k_2 \in \{ 1, \ldots, m \}$,
$J_1$, $J_2 \in (\mathbb{Z}_{\ge 0})^n$ we have
\begin{enumerate}
\item \label{mlhdr:ranking1}
$\partial_j \, u_k > u_k$ and
\item \label{mlhdr:ranking2}
$(u_{k_1})_{J_1} > (u_{k_2})_{J_2}$
implies $\partial_j \, (u_{k_1})_{J_1} > \partial_j \, (u_{k_2})_{J_2}$.
\end{enumerate}
\end{definition}

\begin{remark}
Every ranking $>$ on $R$ is a well-ordering
(cf., e.g., \cite[Ch.~0, Sect.~17, Lemma~15]{Kolchin}),
i.e., every descending sequence of elements of ${\rm Mon}(\Delta) \, u$
terminates.
\end{remark}

\begin{example}\label{mlhdr:ex:degrevlex}
On $K\{ u \}$ (i.e., $m = 1$)
with commuting derivations $\partial_1$, \ldots, $\partial_n$
the \emph{degree-reverse lexicographical ranking} (with
$\partial_1 \, u > \partial_2 \, u > \ldots > \partial_n \, u$)
is defined for $u_{J}$, $u_{J'}$, $J = (j_1, \ldots, j_n)$,
$J' = (j'_1, \ldots, j'_n) \in (\mathbb{Z}_{\ge 0})^n$, by
\[
u_{J} > u_{J'}
\quad :\Longleftrightarrow \quad
\left\{ \begin{array}{l}
j_1 + \ldots + j_n > j'_1 + \ldots + j'_n \quad \mbox{or}\\[1em]
\quad \big(
\phantom{\{} j_1 + \ldots + j_n = j'_1 + \ldots + j'_n \quad \mbox{and} \quad
J \neq J' \quad \mbox{and}\\[1em]
\quad j_i < j'_i \quad \mbox{for}
\quad i = \max \, \{ \, 1 \le k \le n \mid j_k \neq j'_k \, \} \phantom{\}}
\big).
\end{array} \right.
\]
For instance, if $n = 3$, we have
$u_{(1,2,1)} > u_{(1,2,0)} > u_{(2,0,1)}$.
\end{example}

In what follows, we assume that a ranking $>$ on $R = K\{ u_1, \ldots, u_m \}$
is fixed.

\begin{remark}
Let $p_1$, $p_2 \in R$ be two non-constant differential polynomials. If $p_1$
and $p_2$ have the same leader $(u_k)_J$ and the degree of $p_1$ in $(u_k)_J$
is greater than or equal to the degree of $p_2$ in $(u_k)_J$, then the same
pseu\-do-di\-vi\-sion as in (\ref{mlhdr:eq:pseudoreduction}) yields a
remainder which is either zero, or has leader less than $(u_k)_J$, or has
leader $(u_k)_J$ and smaller degree in $(u_k)_J$ than $p_1$.

More generally, if ${\rm ld}(p_1) = \theta \, {\rm ld}(p_2)$ for
some $\theta \in {\rm Mon}(\Delta)$, then this pseu\-do-di\-vi\-sion can be
applied with $p_2$ replaced with $\theta \, p_2$.
Note that, by condition \ref{mlhdr:ranking2} of the definition of a ranking,
we have ${\rm ld}(\theta \, p_2) = \theta \, {\rm ld}(p_2)$, and that, if
$\theta \neq 1$, the degree of $\theta \, p_2$ in $\theta \, {\rm ld}(p_2)$
is one, so that the reduction can be applied without assumption on the
degree of $p_2$ in ${\rm ld}(p_2)$.
Then $c_1$ in (\ref{mlhdr:eq:pseudoreduction}) is again chosen as a
suitable power of ${\rm init}(\theta \, p_2)$. In case $\theta \neq 1$
we have
\[
{\rm init}(\theta \, p_2) = \frac{\partial p_2}{\partial \, {\rm ld}(p_2)}
=: {\rm sep}(p_2),
\]
and this differential polynomial is referred to as the \emph{separant} of $p_2$.

In order not to change the solution set of a differential system,
when $p_1 = 0$ is replaced with $r = 0$, where $r$ is the result of
a reduction of $p_1$ modulo $p_2$ or $\theta \, p_2$ as above, it is
assumed that ${\rm init}(p_2)$ and ${\rm sep}(p_2)$ do not vanish on
the solution set of the system. By definition of the separant and the
discriminant (cf.\ Definition~\ref{mlhdr:de:prep} \ref{mlhdr:discriminant}),
non-vanishing of ${\rm sep}(p_2)$ follows from non-vanishing
of ${\rm disc}(p_2)$, as ensured by the algebraic part of Thomas' algorithm
(cf.\ Remark~\ref{mlhdr:rem:algebraicpart}).
\end{remark}

We assume now that the given differential system is simple as an algebraic
system; it could be one of the systems resulting from the algebraic part of
Thomas' algorithm.

\begin{remark}
The symmetry of the
second derivatives $\partial_i \, \partial_j \, u_k = \partial_j \, \partial_i \, u_k$
(and similarly for higher order derivatives) imposes necessary conditions on
the solvability of a system of partial differential equations. Taking identities
like these into account and forming linear combinations of (derivatives of) the
given equations may produce differential consequences with lower ranked leaders. In
order to obtain a complete set of algebraic conditions on the Taylor coefficients
of an analytic solution, the system has to be augmented by these integrability
conditions in general. If a system of partial differential equations admits a
translation into algebraic conditions on the Taylor coefficients such that no
further integrability conditions have to be taken into account, then
it is said to be \emph{formally integrable}.
\end{remark}

A simple differential system, to be defined in Definition~\ref{mlhdr:de:diffsimple},
will be assumed to be formally integrable. The construction of simple differential
systems, and therefore, the computation of a Thomas decomposition, as presented in
\cite{BaechlerGerdtLangeHegermannRobertz2}, \cite{Robertz6}, employs techniques
which can be traced back to C.~Riquier \cite{Riquier} and M.~Janet \cite{Janet1}.
The main idea is to turn the search for new differential consequences (i.e.,
integrability conditions) into a systematic procedure by singling out for
each differential equation those derivations (called ``non-admissible'' here)
which need to be applied to it in this investigation. The notion of Janet division,
as discussed next, establishes a sense of direction in combining the given
equations and deriving consequences. It is a particular case of
an involutive division on sets of monomials, a concept developed by V.~P.\ Gerdt
and Y.~A.\ Blinkov and others (cf., e.g., \cite{GerdtBlinkov98a}).

\begin{definition}
Given a finite subset $M$ of ${\rm Mon}(\Delta)$,
\emph{Janet division} associates with each $\theta \in M$ a subset of
\emph{admissible derivations} $\mu(\theta, M)$ of $\Delta = \{ \partial_1, \ldots, \partial_n \}$
as follows. Let $\theta = \partial_1^{i_1} \ldots \partial_n^{i_n}$.
Then $\partial_k \in \mu(\theta, M)$ if and only if
\[
i_k = \max \, \{ \, j_k \mid \partial_1^{j_1} \ldots \partial_n^{j_n} \in M \mbox{ with }
j_1 = i_1, \, j_2 = i_2, \ldots, \, j_{k-1} = i_{k-1} \, \}.
\]
The subset $\overline{\mu}(\theta, M) := \Delta \setminus \mu(\theta, M)$
consists of the \emph{non-admissible derivations} for the element $\theta$ of $M$.
\end{definition}

\newlength\mlhdrmultvarwidth
\settowidth{\mlhdrmultvarwidth}{$\delta_{w}$}

\begin{example}\label{mlhdr:ex:janetdivision}
Let $\Delta = \{ \, \partial_1, \partial_2, \partial_3 \, \}$ and
$M = \{ \, \partial_1^2 \, \partial_2, \,
\partial_1^2 \, \partial_3, \,
\partial_2^2 \, \partial_3, \,
\partial_2 \, \partial_3^2 \, \}$.
Then Janet division associates the sets $\mu(\theta, M)$ of
admissible derivations to the elements $\theta \in M$ as indicated
in the following table, where we replace non-admissible derivations
in the set $\Delta$ with the symbol '$*$'.
\[
\begin{array}{rl}
\partial_1^2 \, \partial_2, & \,
\{ \makebox[\mlhdrmultvarwidth][c]{$\partial_1$},
\makebox[\mlhdrmultvarwidth][c]{$\partial_2$},
\makebox[\mlhdrmultvarwidth][c]{$\partial_3$} \}\\[0.5em]
\partial_1^2 \, \partial_3, & \,
\{ \makebox[\mlhdrmultvarwidth][c]{$\partial_1$},
\makebox[\mlhdrmultvarwidth][c]{$*$},
\makebox[\mlhdrmultvarwidth][c]{$\partial_3$} \}\\[0.5em]
\partial_2^2 \, \partial_3, & \,
\{ \makebox[\mlhdrmultvarwidth][c]{$*$},
\makebox[\mlhdrmultvarwidth][c]{$\partial_2$},
\makebox[\mlhdrmultvarwidth][c]{$\partial_3$} \}\\[0.5em]
\partial_2 \, \partial_3^2, & \,
\{ \makebox[\mlhdrmultvarwidth][c]{$*$},
\makebox[\mlhdrmultvarwidth][c]{$*$},
\makebox[\mlhdrmultvarwidth][c]{$\partial_3$} \}
\end{array}
\]
\end{example}

\begin{definition}
A finite subset $M$ of ${\rm Mon}(\Delta)$ is said to be \emph{Janet complete} if
\[
\bigcup_{\theta \in M} {\rm Mon}(\mu(\theta, M)) \, \theta = \bigcup_{\theta \in M} {\rm Mon}(\Delta) \, \theta,
\]
i.e., if every monomial which is divisible by some monomial in $M$ is obtained by
multiplying a certain $\theta \in M$ by admissible derivations for $\theta$ only.
(Recall that the left hand side of the above equation is a disjoint union.)
\end{definition}

\begin{example}
The set $M$ in Example~\ref{mlhdr:ex:janetdivision} is not Janet complete
because, e.g., the monomial $\partial_1 \, \partial_2^2 \, \partial_3$ is
not obtained as a multiple of any $\theta \in M$ when multiplication is
restricted to admissible derivations for $\theta$. By adding this monomial
and the monomial $\partial_1 \, \partial_2 \, \partial_3^2$ to $M$, we
obtain the following Janet complete superset of $M$ in ${\rm Mon}(\Delta)$.
\[
\begin{array}{rl}
\partial_1^2 \, \partial_2, & \,
\{ \makebox[\mlhdrmultvarwidth][c]{$\partial_1$},
\makebox[\mlhdrmultvarwidth][c]{$\partial_2$},
\makebox[\mlhdrmultvarwidth][c]{$\partial_3$} \}\\[0.5em]
\partial_1^2 \, \partial_3, & \,
\{ \makebox[\mlhdrmultvarwidth][c]{$\partial_1$},
\makebox[\mlhdrmultvarwidth][c]{$*$},
\makebox[\mlhdrmultvarwidth][c]{$\partial_3$} \}\\[0.5em]
\partial_1 \, \partial_2^2 \, \partial_3, & \,
\{ \makebox[\mlhdrmultvarwidth][c]{$*$},
\makebox[\mlhdrmultvarwidth][c]{$\partial_2$},
\makebox[\mlhdrmultvarwidth][c]{$\partial_3$} \}\\[0.5em]
\partial_1 \, \partial_2 \, \partial_3^2, & \,
\{ \makebox[\mlhdrmultvarwidth][c]{$*$},
\makebox[\mlhdrmultvarwidth][c]{$*$},
\makebox[\mlhdrmultvarwidth][c]{$\partial_3$} \}\\[0.5em]
\partial_2^2 \, \partial_3, & \,
\{ \makebox[\mlhdrmultvarwidth][c]{$*$},
\makebox[\mlhdrmultvarwidth][c]{$\partial_2$},
\makebox[\mlhdrmultvarwidth][c]{$\partial_3$} \}\\[0.5em]
\partial_2 \, \partial_3^2, & \,
\{ \makebox[\mlhdrmultvarwidth][c]{$*$},
\makebox[\mlhdrmultvarwidth][c]{$*$},
\makebox[\mlhdrmultvarwidth][c]{$\partial_3$} \}
\end{array}
\]
\end{example}

\begin{remark}
Every finite subset $M$ of ${\rm Mon}(\Delta)$ can be augmented to a
Janet complete finite set by adding certain monomials which are products
of some $\theta \in M$ and a monomial which is divisible by at least one
non-admissible derivation for $\theta$.

For more details on Janet division, we refer to, e.g.,
\cite{GerdtBlinkov98a}, \cite{BaechlerGerdtLangeHegermannRobertz2}, \cite{Robertz6}.
\end{remark}

Each equation $p_i = 0$ in a differential system is assigned the set of
admissible derivations $\mu(\theta_i, M_k)$,
where ${\rm ld}(p_i) = \theta_i \, u_k$ and
\begin{equation}\label{mlhdr:eq:mk}
M_k := \{ \, \theta \in {\rm Mon}(\Delta) \mid \theta \, u_k \in \{ \, {\rm ld}(p_1), \ldots, {\rm ld}(p_s) \, \} \, \}
\end{equation}
is the set of all monomials which define leaders of the
equations $p_1 = 0$, \ldots, $p_s = 0$ in the system
involving the same differential indeterminate $u_k$.
We refer to $d \, p_i$ for $d \in {\rm Mon}(\mu(\theta_i, M_k))$
as the \emph{admissible derivatives} of $p_i$.

\medskip

Formal integrability of a differential system is then decided by applying to
each equation $p_i = 0$ every of its non-admissible derivations
$d \in \overline{\mu}(\theta_i, M_k)$ and computing the
pseu\-do-re\-main\-der of $d \, p_i$ modulo $p_1$, \ldots, $p_s$ and their
admissible derivatives. The restriction of the pseu\-do-di\-vi\-sion to
admissible derivatives requires $M_k$ to be Janet complete.
If one of these pseu\-do-re\-main\-ders is non-zero, then it is added as a
new equation to the system, and the augmented system has to be treated by
the algebraic part of Thomas' algorithm again.

\begin{definition}
\renewcommand{\theenumi}{\alph{enumi})}
\renewcommand{\labelenumi}{\theenumi}
\renewcommand{\theenumii}{\roman{enumii}}
A system of partial differential equations $\{ \, p_1 = 0, \ldots, p_s = 0 \, \}$,
where $p_1$, \ldots, $p_s \in R \setminus K$, is said to be \emph{passive} if the
following two conditions hold for ${\rm ld}(p_1) = \theta_1 \, u_{k_1}$, \ldots,
${\rm ld}(p_s) = \theta_s \, u_{k_s}$, where
$\theta_i \in {\rm Mon}(\Delta)$, $k_i \in \{ 1, \ldots, m \}$.
\begin{enumerate}
\item \label{mlhdr:passive1}
For all $k \in \{ 1, \ldots, m \}$, the set $M_k$ defined in (\ref{mlhdr:eq:mk})
is Janet complete.
\item \label{mlhdr:passive2}
For all $i \in \{ 1, \ldots, s \}$ and all $d \in \overline{\mu}(\theta_i, M_{k_i})$,
the pseu\-do-re\-main\-der of $d \, p_i$ modulo $p_1$, \ldots, $p_s$ and
their admissible derivatives is zero.
\end{enumerate}
\end{definition}

\begin{definition}\label{mlhdr:de:diffsimple}
\renewcommand{\theenumi}{\alph{enumi})}
\renewcommand{\labelenumi}{\theenumi}
\renewcommand{\theenumii}{\roman{enumii}}
A differential system $S$, defined over $R$, as in (\ref{mlhdr:eq:diffsys})
is said to be \emph{simple}
(with respect to $>$) if the following three conditions hold.
\begin{enumerate}
\item \label{mlhdr:diffsimple1}
The system $S$ is simple as an algebraic system (in the
finitely many variables $(u_i)_J$ which occur in the equations
and inequations of $S$, totally ordered by $>$).
\item \label{mlhdr:diffsimple2}
The system $\{ \, p_1 = 0, \ldots, p_s = 0 \, \}$ is passive.
\item \label{mlhdr:diffsimple3}
The left hand sides of the inequations $q_1 \neq 0$, \ldots, $q_t \neq 0$
equal their pseu\-do-re\-main\-ders modulo $p_1$, \ldots, $p_s$
and their derivatives.
\end{enumerate}
\end{definition}

\begin{definition}
Let $S$ be a differential system, defined over $R$.
A \emph{Thomas decomposition} of $S$ (or of ${\rm Sol}_{\Omega}(S)$)
with respect to $>$
is a collection of finitely many simple differential systems
$S_1$, \ldots, $S_r$, defined over $R$, such that
${\rm Sol}_{\Omega}(S)$ is the disjoint union of the solution sets
${\rm Sol}_{\Omega}(S_1)$, \ldots, ${\rm Sol}_{\Omega}(S_r)$.
\end{definition}

\begin{remark}\label{mlhdr:rem:differentialpart}
Given $S$ as in (\ref{mlhdr:eq:diffsys}) and a ranking on $R$, a
Thomas decomposition of $S$ with respect to $>$ can be computed
by interweaving the algebraic part discussed in Subsection~\ref{mlhdr:subsec:algsystem}
and differential reduction and completion with respect to Janet division.

First of all, a Thomas decomposition of $S$, considered as an
algebraic system, is computed. Each of the resulting simple algebraic
systems is then treated as follows. Differential pseu\-do-di\-vi\-sion
is applied to pairs of distinct equations with leaders $\theta_1 \, u_k$
and $\theta_2 \, u_k$ such that $\theta_1 \mid \theta_2$ until either
a non-zero pseu\-do-re\-main\-der is obtained or no such further
reductions are possible. Non-zero pseu\-do-re\-main\-ders are added to
the system, the algebraic part of Thomas' algorithm is applied again,
and the process is repeated. Once the system is auto-reduced in this
sense, then it is possibly augmented with certain derivatives of equations
so that the sets $M_k$ defined in (\ref{mlhdr:eq:mk}) are Janet complete.
Then it is checked whether the system is passive. If a non-zero remainder
is obtained by a pseu\-do-di\-vi\-sion of a non-admissible derivative
modulo the equations and their admissible derivatives, then the
algebraic part of Thomas' algorithm is applied again to the augmented
system. Otherwise, the system is passive. Finally, the left hand side
of each inequation is replaced with its pseu\-do-re\-main\-der modulo
the equations and their derivatives, in order to ensure
condition~\ref{mlhdr:diffsimple3} of Definition~\ref{mlhdr:de:diffsimple}.
The main reason why this procedure terminates is Dickson's Lemma, which
shows that the ascending sequence of ideals of the
semigroup ${\rm Mon}(\Delta)$ formed by the monomials $\theta$ defining
leaders of equations (for each differential indeterminate) becomes
stationary after finitely many steps.

For more details on the differential part of Thomas' algorithm,
we refer to \cite{BaechlerGerdtLangeHegermannRobertz2}, \cite{LangeHegermannThesis},
and \cite[Subsect.~2.2.2]{Robertz6}.

An implementation of Thomas' algorithm for differential systems has been
developed by M.~Lange-Hegermann as Maple package {\tt DifferentialThomas}
\cite{ThomasPackages}.
\end{remark}

We also use a simpler notation for the
indeterminates $(u_k)_J$ of the differential polynomial ring.
In case $m = 1$ we use the symbol $u$ as a synonym for $u_1$. In addition,
if the derivations $\partial_1$, $\partial_2$, $\partial_3$ represent
the partial differential operators with respect to $x$, $y$, $z$,
respectively, then we write
\[
u_{\underbrace{\mbox{\scriptsize $x, \ldots, x$}}_{i},
\underbrace{\mbox{\scriptsize $y, \ldots, y$}}_{j},
\underbrace{\mbox{\scriptsize $z, \ldots, z$}}_{k}}
\]
instead of $u_{(i,j,k)}$.

\medskip

When displaying a simple differential system we indicate next to each
equation its set of admissible derivations.

\begin{example}
Let us consider the ordinary differential equation (which is
discussed in \cite[Example in Sect.~4.7]{Ince})
\[
\left( \frac{\partial u}{\partial x} \right)^3 - 4 \, x \, u(x) \, \frac{\partial u}{\partial x} + 8 \, u(x)^2 = 0.
\]
The left hand side is represented by the
element $p := u_x^3 - 4 \, x \, u \, u_x + 8 \, u^2$ of the
differential polynomial ring $R = K\{ u \}$ with one derivation $\partial_x$,
where $K = \mathbb{Q}(x)$ is the field of rational functions in $x$,
endowed with differentiation with respect to $x$.

The initial of $p$ is constant, the separant of $p$ is $3 \, u_x^2 - 4 \, x \, u$.
The algebraic part of Thomas' algorithm only distinguishes the
cases whether the discriminant of $p$ vanishes or not. We have
\[
{\rm disc}(p) = -{\rm res}(p, {\rm sep}(p), u_x) = -64 \, \underline{u}^3 \, (27 \, \underline{u} - 4 \, x^3).
\]
This case distinction leads to the Thomas decomposition
\begin{center}
\fbox{\;
$\begin{array}{rcll}
\underline{u_x}^3 - 4 \, x \, u \, \underline{u_x} + 8 \, u^2 & = & 0, & \,
\{ \makebox[\mlhdrmultvarwidth][c]{$\partial_x$} \}\\[1em]
(27 \, \underline{u} - 4 \, x^3) \, \underline{u} & \neq & 0 &
\end{array}$
\;}
\ 
\fbox{\;
$\begin{array}{rcll}
\phantom{x}\\[1em]
(27 \, \underline{u} - 4 \, x^3) \, \underline{u} & = & 0, & \,
\{ \makebox[\mlhdrmultvarwidth][c]{$\partial_x$} \}
\end{array}$
\;}
\end{center}
Since both systems contain only one equation, no differential reductions
are necessary. The second simple system could be split into two with
equations $27 \, u - 4 \, x^3 = 0$ and $u = 0$, respectively.
The solutions of the first simple system are given by
$u(x) = c \, (x - c)^2$, where $c$ is an arbitrary non-zero constant.
The solutions $u(x) = 0$ and $u(x) = \frac{4}{27} \, x^3$ of
the second simple system are called \emph{singular solutions},
the latter one being an envelope of the general solution.
\end{example}

\begin{example}
Let us compute a Thomas decomposition of the system of (nonlinear)
partial differential equations
\[
\left\{
\begin{array}{rcl}
\displaystyle
\frac{\partial^2 u}{\partial x^2} - \frac{\partial^2 u}{\partial y^2} & = & 0,\\[1em]
\displaystyle
\frac{\partial u}{\partial x} - u^2 & = & 0
\end{array} \right.
\]
for one unknown function $u(x, y)$. The left hand sides are
expressed as elements $p_1 := u_{x,x} - u_{y,y}$ and $p_2 := u_x - u^2$
of the differential polynomial ring
$R = \mathbb{Q}\{ u \}$ with commuting derivations $\partial_x$, $\partial_y$.
We choose the degree-reverse lexicographical ranking $>$ on $R$ with
$\partial_x \, u > \partial_y \, u$ (cf.\ Example~\ref{mlhdr:ex:degrevlex}).

Since the monomial $\partial_x$ defining the leader of $p_2$ divides
the monomial $\partial_x^2$ defining the leader of $p_1$, differential
pseu\-do-di\-vi\-sion is applied and $p_1$ is replaced with
\[
p_3 := p_1 - \partial_x \, p_2 - 2 \, u \, p_2 = -u_{y,y} + 2 \, u^3.
\]
Janet division associates the sets of admissible derivations to the
equations of the resulting system as follows:
\[
\left\{
\begin{array}{rcll}
\underline{u_x} - u^2 & = & 0, & \,
\{ \makebox[\mlhdrmultvarwidth][c]{$\partial_x$},
\makebox[\mlhdrmultvarwidth][c]{$\partial_y$} \}\\[1em]
\underline{u_{y,y}} - 2 \, u^3 & = & 0, & \,
\{ \makebox[\mlhdrmultvarwidth][c]{$*$},
\makebox[\mlhdrmultvarwidth][c]{$\partial_y$} \}
\end{array} \right.
\]
The set of monomials $\{ \, \partial_x, \partial_y^2 \, \}$ defining
the leaders $u_x$ and $u_{y,y}$ is Janet complete. The check whether
the above system is passive involves the following reduction:
\[
\partial_x \, p_3 + \partial_y^2 \, p_2 - 6 \, u^2 \, p_2 - 2 \, u \, p_3 =
-2 \, (\underline{u_y} + u^2) \, (\underline{u_y} - u^2).
\]
This non-zero remainder is a differential consequence which is added
as an equation to the system. In fact, the system can be split into
two systems according to the given factorization. For both systems a
differential reduction of $p_3$ modulo the chosen factor is applied
because the monomial $\partial_y$ defining the new leader divides the
monomial $\partial_{y,y}$ defining ${\rm ld}(p_3)$. In both cases the
remainder is zero, the sets of monomials defining leaders are Janet
complete, and the passivity check confirms formal integrability. We
obtain the Thomas decomposition
\begin{center}
\fbox{\;
$\begin{array}{rcll}
\underline{u_x} - u^2 & = & 0, & \,
\{ \makebox[\mlhdrmultvarwidth][c]{$\partial_x$},
\makebox[\mlhdrmultvarwidth][c]{$\partial_y$} \}\\[1em]
\underline{u_y} + u^2 & = & 0, & \,
\{ \makebox[\mlhdrmultvarwidth][c]{$*$},
\makebox[\mlhdrmultvarwidth][c]{$\partial_y$} \}\\[1em]
\phantom{x}
\end{array}$
\;}
\qquad
\fbox{\;
$\begin{array}{rcll}
\underline{u_x} - u^2 & = & 0, & \,
\{ \makebox[\mlhdrmultvarwidth][c]{$\partial_x$},
\makebox[\mlhdrmultvarwidth][c]{$\partial_y$} \}\\[1em]
\underline{u_y} - u^2 & = & 0, & \,
\{ \makebox[\mlhdrmultvarwidth][c]{$*$},
\makebox[\mlhdrmultvarwidth][c]{$\partial_y$} \}\\[1em]
\underline{u} & \neq & 0. &
\end{array}$
\;}
\end{center}
If the above factorization is ignored, then the
discriminant of $p_4 := u_y^2 - u^4$ needs to be considered, which
implies vanishing or non-vanishing of the separant $2 \, u_y$.
This case distinction leads to a different Thomas decomposition.
\end{example}

A Thomas decomposition of a differential system is not uniquely
determined, as the previous example shows (cf.\ also
Remark~\ref{mlhdr:rem:thomasdecompnonunique} for the algebraic case).
In the special case of a system $S$ of \emph{linear} partial
differential equations no case distinctions are necessary, and
the single simple system in any Thomas decomposition of $S$ is
a Janet basis for $S$ (cf., e.g., \cite{Janet1}, \cite{Pommaret1},
\cite{GerdtBlinkov98a}, \cite{Robertz6}). Pseu\-do-re\-duc\-tion
of a differential polynomial
modulo the equations of a simple differential system and their
derivatives decides membership to the corresponding saturation
ideal (cf.\ also Proposition~\ref{mlhdr:prop:algebraicmembership}).

\begin{proposition}[\cite{Robertz6}, Prop.~2.2.50]\label{mlhdr:prop:differentialmembership}
Let $S$ be a simple differential system, defined over $R$, with
equations $p_1 = 0$, \ldots, $p_s = 0$. Moreover, let $E$ be the
differential ideal of $R$ generated by $p_1$, \ldots, $p_s$ and
define the product $q$ of the initials and separants of all $p_1$,
\ldots, $p_s$. Then $E : q^{\infty}$ is a radical differential ideal.
Given $p \in R$, we have $p \in E : q^{\infty}$ if and only if the
pseu\-do-re\-main\-der of $p$ modulo $p_1$, \ldots, $p_s$
and their derivatives is zero.
\end{proposition}

Similarly to the algebraic case, the Nullstellensatz for analytic functions
(due to J.~F.\ Ritt and H.~W.\ Raudenbush, cf.\ \cite[Sects.\ II.7--11, IX.27]{Ritt})
establishes a one-to-one correspondence of solution sets $V := {\rm Sol}_{\Omega}(S)$
of systems of partial differential equations $S = \{ \, p_1 = 0, \ldots, p_s = 0 \, \}$
for $m$ unknown functions, defined over $R$, and their vanishing ideals
in $R = K\{ u_1, \ldots, u_m \}$
\[
\mathcal{I}_R(V) := \{ \, p \in R \mid p(f) = 0 \mbox{ for all } f \in V \, \}.
\]
These are the radical differential ideals of $R$. The Nullstellensatz
implies that, with the notation of Proposition~\ref{mlhdr:prop:differentialmembership},
we have $\mathcal{I}_R({\rm Sol}_{\Omega}(S)) = E : q^{\infty}$.

\medskip

The following proposition allows to decide whether a given differential
equation $p = 0$ is a consequence of a (not necessarily simple) differential
system $S$ by applying pseu\-do-di\-vi\-sion to $p$ modulo each of the simple
systems in a Thomas decomposition of $S$.
It follows from the previous proposition and the Nullstellensatz and it
also applies to algebraic systems by ignoring the separants.

\begin{proposition}[\cite{Robertz6}, Prop.~2.2.72]\label{mlhdr:prop:intersection}
Let $S$ be a (not necessarily simple) differential system
as in (\ref{mlhdr:eq:diffsys})
and $S_1$, \ldots, $S_r$ a Thomas decomposition of $S$ with respect to
any ranking on $R$. Moreover, let $E$ be the differential ideal of $R$
generated by $p_1$, \ldots, $p_s$ and define the product $q$ of $q_1$, \ldots, $q_t$.
For $i \in \{ 1, \ldots, r \}$, let $E^{(i)}$ be the differential ideal
of $R$ generated by the equations in $S_i$ and define the product $q^{(i)}$
of the initials and separants of all these equations. Then we have
\[
\sqrt{E : q^{\infty}} = \left( E^{(1)} : (q^{(1)})^{\infty} \right) \cap \ldots \cap
\left( E^{(r)} : (q^{(r)})^{\infty} \right).
\]
\end{proposition}

\section{Elimination}
\label{mlhdr:sec:elimination}

Thomas' algorithm can be used to solve various differential elimination problems.
This section presents results on certain rankings on the differential polynomial
ring $R = K\{ u_1, \ldots, u_m \}$ which allow to compute all differential
consequences of a given differential system involving only a specified subset of
the differential indeterminates $u_1$, \ldots, $u_m$. In other words, this technique
allows to determine all differential equations which are satisfied by certain
components of the solution tuples. We adopt the notation from the previous section.

\begin{definition}
Let $I_1$, $I_2$, \ldots, $I_k$
form a partition of $\{ 1, 2, \ldots, m \}$ such that
$i_1 \in I_{j_1}$, $i_2 \in I_{j_2}$, $i_1 \le i_2$ implies
$j_1 \le j_2$.
Let $B_j := \{ u_i \mid i \in I_j \}$, $j = 1$, \ldots, $k$.
Moreover, fix some degree-reverse lexicographical ordering $>$
on ${\rm Mon}(\Delta)$.
Then the \emph{block ranking} on $R$ \emph{with blocks $B_1$, \ldots, $B_k$}
(with $u_1 > u_2 > \ldots > u_m$) is defined for
$\theta_1 \, u_{i_1}$, $\theta_2 \, u_{i_2} \in {\rm Mon}(\Delta) \, u$,
where $u_{i_1} \in B_{j_1}$, $u_{i_2} \in B_{j_2}$, by
\[
\theta_1 \, u_{i_1} > \theta_2 \, u_{i_2}
\quad :\Longleftrightarrow \quad
\left\{ \begin{array}{l}
j_1 < j_2 \quad \mbox{or} \quad \Big( \phantom{\{} j_1 = j_2 \quad
\mbox{and} \quad \big(\, \theta_1 > \theta_2 \quad \mbox{or}\\[1em]
\quad (\, \theta_1 = \theta_2 \quad \mbox{and} \quad i_1 < i_2
\, ) \, \big) \phantom{\}} \Big).
\end{array} \right.
\]
Such a ranking is said to satisfy $B_1 \gg B_2 \gg \ldots \gg B_k$.
\end{definition}

\begin{example}
With respect to the block ranking on $K\{ u_1, u_2, u_3 \}$ with
blocks $\{ u_1 \}$, $\{ u_2, u_3 \}$ (and $u_1 > u_2 > u_3$) we have
$(u_1)_{(0,1)} > u_1 > (u_2)_{(1,2)} > (u_3)_{(1,2)} > (u_2)_{(0,1)}$.
\end{example}

In the situation of the previous definition,
for every $i \in \{ 1, \ldots, k \}$, we consider
$K\{ B_i, \ldots, B_k \} := K\{ u \mid u \in B_i \cup \ldots \cup B_k \}$
as a differential subring of $R$, endowed with the restrictions of
the derivations $\partial_1$, \ldots, $\partial_n$ to $K\{ B_i, \ldots, B_k \}$.

\medskip

For any algebraic or differential system $S$ we denote by $S^{=}$ (resp.\ $S^{\neq}$)
the set of the left hand sides of all equations (resp.\ inequations) in $S$.

\begin{proposition}[\cite{Robertz6}, Prop.~3.1.36]\label{mlhdr:prop:elimsimple}
Let $S$ be a simple differential system, defined over $R$, with respect to
a block ranking with blocks $B_1$, \ldots, $B_k$. Moreover, let $E$ be the
differential ideal of $R$ generated by $S^{=}$ and $q$ the product of the
initials and separants of all elements of $S^{=}$. For every $i \in \{ 1, \ldots, k \}$,
let $E_i$ be the differential ideal of $K\{ B_i, \ldots, B_k \}$ generated
by $P_i := S^{=} \cap K\{ B_i, \ldots, B_k \}$ and let $q_i$ be the product
of the initials and separants of all elements of $P_i$. Then, for every
$i \in \{ 1, \ldots, k \}$, we have
\[
(E : q^{\infty}) \cap K\{ B_i, \ldots, B_k \} = E_i : q_i^{\infty}\,.
\]
\end{proposition}

In other words, the differential equations implied by $S$ which involve only
the differential indeterminates in $B_i \cup \ldots \cup B_k$ are precisely
those whose pseu\-do-re\-main\-ders modulo the elements
of $S^{=} \cap K\{ B_i, \ldots, B_k \}$ and their derivatives are zero.

\begin{example}
The Cauchy-Riemann equations for a complex function of $z = x+i \, y$ with real part
$u$ and imaginary part $v$ are
\[
\left\{
\begin{array}{rcl}
\displaystyle
\frac{\partial u}{\partial x} - \frac{\partial v}{\partial y} & = & 0,\\[1em]
\displaystyle
\frac{\partial u}{\partial y} + \frac{\partial v}{\partial x} & = & 0.
\end{array}
\right.
\]
The left hand sides are represented by the elements
$p_1 := u_x - v_y$ and $p_2 := u_y + v_x$ of the differential polynomial
ring $R = \mathbb{Q}\{ u, v \}$ with derivations $\partial_x$ and $\partial_y$.
Choosing a block ranking on $R$ satisfying $\{ u \} \gg \{ v \}$, the
passivity check yields the equation
\[
\partial_x \, p_2 - \partial_y \, p_1 = v_{x,x} + v_{y,y} = 0.
\]
Similarly, the choice of a block ranking on $R$ satisfying $\{ v \} \gg \{ u \}$
yields the consequence $u_{x,x} + u_{y,y} = 0$. These computations confirm that
the real and imaginary parts of a holomorphic function are harmonic functions.
\end{example}

\begin{corollary}[\cite{Robertz6}, Cor.~3.1.37]\label{mlhdr:cor:elimination}
Let $S$ be a (not necessarily simple) differential system, defined over $R$,
and $S_1$, \ldots, $S_r$ a Thomas decomposition of $S$ with respect to
a block ranking with blocks $B_1$, \ldots, $B_k$. Moreover, let $E$ be the
differential ideal of $R$ generated by $S^{=}$ and $q$ the product of all
elements of $S^{\neq}$. Let $i \in \{ 1, \ldots, k \}$ be fixed.
For every $j \in \{ 1, \ldots, r \}$, let $E^{(j)}$ be the
differential ideal of $K\{ B_i, \ldots, B_k \}$ generated by
$P_j := S_j^{=} \cap K\{ B_i, \ldots, B_k \}$ and let $q^{(j)}$ be the
product of the initials and separants of all elements of $P_j$. Then we have
\[
\sqrt{E : q^{\infty}} \, \cap \, K\{ B_i, \ldots, B_k \} =
(E_1 : q_1^{\infty}) \, \cap \, \ldots \, \cap \, (E_r : q_r^{\infty})\,.
\]
\end{corollary}

\section{Control-theoretic applications}
\label{mlhdr:sec:control}

In order to apply the Thomas decomposition method to nonlinear control
systems, we assume that the control system is given by differential
equations and inequations whose left hand sides are polynomials.
Structural information about certain configurations of the control system
is obtained from each simple system of a Thomas decomposition of the
given differential equations and inequations. The choice of ranking on the
differential polynomial ring depends on the question at hand, although a
Thomas decomposition with respect to any ranking, e.g., the degree-reverse
lexicographical ranking, may give hints on how to adapt the ranking for
further investigations in a certain direction.

Let $R = K\{ U \}$ be the differential polynomial ring in the
differential indeterminates $U := \{ u_1, \ldots, u_m \}$ over a
differential field $K$ of (complex) meromorphic functions on an open and
connected subset $\Omega$ of $\mathbb{C}^n$ (cf.\ Subsection~\ref{mlhdr:subsec:diffsystem}).
(No distinction is made a priori between state variables, input, output, etc.)

We assume that $S$ is a simple differential system, defined over $R$,
with respect to some ranking $>$.
Let $E$ be the differential ideal of $R$ generated by the set $S^{=}$ of
the left hand sides of the equations in $S$ and define the product $q$ of
the initials and separants of all elements of $S^{=}$.

\begin{definition}\label{mlhdr:de:observable}
Let $x \in U$ and $Y \subseteq U \setminus \{ x \}$.
Then $x$ is said to be \emph{observable with respect to $Y$} if there
exists $p \in (E : q^{\infty}) \setminus \{ 0 \}$ such that $p$ is a
polynomial in $x$ (not involving any proper derivative of $x$) with
coefficients in $K\{ Y \}$ and such that neither its leading coefficient
nor $\partial p / \partial x$ is an element of $E : q^{\infty}$.
\end{definition}

\begin{remark}
Let $p$ be a polynomial as in the previous definition. Then the
implicit function theorem allows to solve $p = 0$ locally for $x$ in the
sense that the component of $(f_1, \ldots, f_m) \in {\rm Sol}_{\Omega}(S)$
corresponding to $x$ can locally be expressed as an analytic function of
the components corresponding to the differential indeterminates in $Y$.

If $>$ satisfies $U \setminus (Y \cup \{ x \}) \gg \{ x \} \gg \{ Y \}$, then
by Proposition~\ref{mlhdr:prop:elimsimple}, there exists a polynomial $p$
in $(E : q^{\infty}) \setminus \{ 0 \}$ as above if and only if there
exists such a polynomial in $S^{=} \cap K\{ Y \cup \{ x \} \}$.
For a not necessarily simple differential system $S$, a Thomas decomposition
with respect to a ranking as above allows to decide the existence of such a
polynomial among the left hand sides of the differential consequences of $S$
by inspecting each simple system (cf.\ Corollary~\ref{mlhdr:cor:elimination}).
\end{remark}

\begin{definition}\label{mlhdr:de:flatoutput}
A subset $Y$ of $U$ is called a \emph{flat output} of $S$
if $(E : q^{\infty}) \cap K\{ Y \} = \{ 0 \}$ and every $x \in U \setminus Y$
is observable with respect to $Y$.
\end{definition}

\begin{remark}
Let $>$ satisfy $U \setminus Y \gg Y$. Then Proposition~\ref{mlhdr:prop:elimsimple}
allows to decide whether the conditions in Definition~\ref{mlhdr:de:flatoutput}
are satisfied by checking that $S^{=} \cap K \{ Y \} = \emptyset$ holds and
that for every $x \in U \setminus Y$ there exists a polynomial
$p \in S^{=} \cap K\{ Y \cup \{ x \} \}$ satisfying the conditions in
Definition~\ref{mlhdr:de:observable}.

If the differential ideal $I := E : q^{\infty}$ is prime, then the field of
fractions ${\rm Quot}(R / I)$ can be considered as a differential extension
field of $K$. Let us assume that $Y$ is a flat output of $S$ and let $L$ be
the differential subfield of ${\rm Quot}(R / I)$ which is generated by
$\{ \, y + I \mid y \in Y \, \}$. Then, by Definition~\ref{mlhdr:de:flatoutput},
$L / K$ is a purely differentially transcendental extension of differential
fields, and for every $x \in U \setminus Y$, the element $x + I$ of
${\rm Quot}(R / I)$ is algebraic over $L$. Hence, $\{ \, y + I \mid y \in Y \, \}$
is a differential transcendence basis of ${\rm Quot}(R / I) / K$,
and the system is flat in the sense of \cite[Sect.~3.2]{FliessLevineMartinRouchon}.
\end{remark}

\begin{remark}
Following \cite{FliessLevineMartinRouchon}, a system which is defined by a
differential field extension is called flat if it is equivalent by endogenous
feedback to a system which is defined by a purely differentially transcendental
extension of differential fields. As opposed to checking whether $Y$ is a flat
output of $S$ using the method described above, deciding whether $S$ is flat is
a difficult problem in general.
\end{remark}

\medskip

As a first illustration of how differential elimination methods can be applied
to nonlinear control systems, we consider inversion, i.e., the problem of
expressing the input variables in terms of the output variables (and their
derivatives).

\begin{remark}
Using the same notation as above, we assume that disjoint subsets $Y$ and $Z$
of $U$ are specified, where the differential indeterminates in $Y$ and $Z$ are
interpreted as input and output variables of the system, respectively. We
achieve an \emph{inversion} of the system $S$ if and only if we can exhibit for
each $z \in Z$ a $p \in (E : q^{\infty}) \setminus \{ 0 \}$ such
that $p$ is a polynomial in $z$ (not involving any proper derivative of $z$) with
coefficients in $K\{ Y \}$ and such that neither its leading coefficient
nor $\partial p / \partial z$ is an element of $E : q^{\infty}$.

If $>$ satisfies $U \setminus (Y \cup \{ z \}) \gg \{ z \} \gg Y$, then
by Proposition~\ref{mlhdr:prop:elimsimple}, there exists such a polynomial $p$
in $(E : q^{\infty}) \setminus \{ 0 \}$ if and only if there
exists such a polynomial in $S^{=} \cap K\{ Y \cup \{ z \} \}$.
A block ranking $U \setminus (Y \cup Z) \gg Z \gg Y$ may allow to find
such polynomials $p$ for all $z \in Z$ by computing only one Thomas decomposition
(cf.\ the following example),
for instance, if all these polynomials $p$ have degree one.
\end{remark}

\medskip

For displaying simple differential systems resulting from Thomas
decompositions in a concise way, we use the following
command {\tt Print}, which makes use of both the Maple packages
{\tt Janet} \cite{BlinkovCidetal} and {\tt DifferentialThomas}
\cite{ThomasPackages}, where {\tt ivar} and {\tt dvar} are the
lists of independent and dependent variables, respectively.

\medskip

\begin{maplegroup}
\begin{mapleinput}
\mapleinline{active}{1d}{with(Janet):}{}
\end{mapleinput}
\end{maplegroup}
\begin{maplegroup}
\begin{mapleinput}
\mapleinline{active}{1d}{Print := S->Diff2Ind(
PrettyPrintDifferentialSystem(S), ivar, dvar):}{}
\end{mapleinput}
\end{maplegroup}

\medskip

The sets of admissible derivations for the equations in a simple system
are not reproduced here. Note that the implementation uses factorization
and may, for convenience, return simple systems containing several
inequations with the same leader (thus, not strictly complying with
condition~\ref{mlhdr:algsimple2} of Definition~\ref{mlhdr:de:algsimple}).

\begin{example}
The following system of ordinary differential equations models a
unicycle as described in \cite[Examples~3.20, 4.18, 5.10]{ConteMoogPerdon}
(cf.\ also, e.g., \cite[Example~2.35]{NijmeijerVanDerSchaft}).
\[
\left\{
\begin{array}{rcl}
\dot{x}_1 & = & \cos(x_3) \, u_1,\\[0.5em]
\dot{x}_2 & = & \sin(x_3) \, u_1,\\[0.5em]
\dot{x}_3 & = & u_2.
\end{array}
\right.
\]
Here $x_1$, $x_2$, $x_3$ are considered as state variables, where $(x_1, x_2)$
is the position of the middle of the axis in the plane and $x_3$ the angle
of its rotation, and the velocities $u_1$, $u_2$ are considered as inputs.
Moreover, the following outputs $y_1$, $y_2$ are given:
\[
\left\{
\begin{array}{rcl}
y_1 & = & x_1\,,\\[0.5em]
y_2 & = & x_2\,.
\end{array}
\right.
\]
The task is to try to invert the system, i.e., to express $u_1$, $u_2$ in
terms of $y_1$, $y_2$ and their derivatives.

\medskip

In order to translate the given equations into differential polynomials,
we represent $\cos(x_3)$ and $\sin(x_3)$ by differential
indeterminates ${\it cx3}$ and ${\it sx3}$ and add the
generating relations
\[
{\it cx3}^2 + {\it sx3}^2 = 1, \quad
{\it cx3}_t = -{\it sx3} \, (x_3)_t, \quad
{\it sx3}_t = {\it cx3} \, (x_3)_t
\]
to the system. More precisely speaking, we adjoin to the
differential polynomial ring $\mathbb{Q}\{ x_1, x_2, x_3, u_1, u_2, y_1, y_2 \}$
with derivation $\partial_t$ the
differential indeterminates ${\it cx3}$ and ${\it sx3}$ and
define the differential ideal $E$ of the resulting differential polynomial ring
which is generated by $(x_1)_t - {\it cx3} \, u_1$,
$(x_2)_t - {\it sx3} \, u_1$, $(x_3)_t - u_2$,
${\it cx3}^2 + {\it sx3}^2 - 1$,
${\it cx3}_t + {\it sx3} \, (x_3)_t$
and ${\it sx3}_t - {\it cx3} \, (x_3)_t$.
We then apply elimination properties of the
differential Thomas decomposition method to $\sqrt{E}$
(see Proposition~\ref{mlhdr:prop:intersection}
and Corollary~\ref{mlhdr:cor:elimination}).

\medskip

(Alternatively, if one accepts neglecting
the particular case of movement of the unicycle in the direction
of the $x_2$-coordinate axis without rotation, one could assume that $\cos(x_3)$
is not the zero function, multiply both sides of the equation $\dot{x}_3 = u_2$
by $\cos(x_3)$,
read the resulting left hand side, using the chain rule,
as the derivative of $\sin(x_3)$,
and obtain the equation ${\it sx3}_t = {\it cx3} \, u_2$.
This would allow to dispose of the differential indeterminate $x_3$
and the computation of the Thomas decomposition below would essentially
yield the first five of the seven simple systems below.)

\medskip

In the following computations concerning the model of a unicycle
the differential polynomial ring
is $\mathbb{Q}\{ x_1, x_2, {\it cx3}, {\it sx3}, x_3, u_1, u_2, y_1, y_2 \}$
with one derivation $\partial_t$.

\bigskip

\begin{maplegroup}
\begin{mapleinput}
\mapleinline{active}{1d}{with(DifferentialThomas):}{}
\end{mapleinput}
\end{maplegroup}
\begin{maplegroup}
\begin{mapleinput}
\mapleinline{active}{1d}{ivar := [t]:}{}
\end{mapleinput}
\end{maplegroup}
\begin{maplegroup}
\begin{mapleinput}
\mapleinline{active}{1d}{dvar := [x1,x2,cx3,sx3,x3,u1,u2,y1,y2]:}{}
\end{mapleinput}
\end{maplegroup}
\medskip
\noindent
We specify the block ranking $>$ satisfying $\{ x_1, x_2, {\it cx3}, {\it sx3}, x_3 \} \gg \{ u_1, u_2 \} \gg \{ y_1, y_2 \}$
as well as $x_1 > x_2 > {\it cx3} > {\it sx3} > x_3$
and $u_1 > u_2$ and $y_1 > y_2$.
\medskip

\begin{maplegroup}
\begin{mapleinput}
\mapleinline{active}{1d}{ComputeRanking(ivar,
[[x1,x2,cx3,sx3,x3],[u1,u2],[y1,y2]]):}{}
\end{mapleinput}
\end{maplegroup}
\medskip
\noindent
If the left hand sides of the system are written in jet notation,
then a conversion into the format expected by the package
{\tt DifferentialThomas} is accomplished by the following sequence
of commands.
\medskip

\begin{maplegroup}
\begin{mapleinput}
\mapleinline{active}{1d}{L := [x1[t]-cx3*u1, x2[t]-sx3*u1, x3[t]-u2,
y1-x1, y2-x2, cx3\symbol{94}2+sx3\symbol{94}2-1, cx3[t]+sx3*x3[t],
sx3[t]-cx3*x3[t]];
}{}
\end{mapleinput}
\end{maplegroup}
\begin{maplegroup}
\begin{mapleinput}
\mapleinline{active}{1d}{LL := Diff2JetList(Ind2Diff(L, ivar, dvar));
}{}
\end{mapleinput}
\mapleresult
\begin{maplelatex}
\mapleinline{inert}{2d}{LL := [x1[1]-cx3[0]*u1[0], x2[1]-sx3[0]*u1[0], x3[1]-u2[0], y1[0]-x1[0], y2[0]-x2[0], cx3[0]^2+sx3[0]^2-1, cx3[1]+sx3[0]*x3[1], sx3[1]-cx3[0]*x3[1]]}{
\[
\begin{array}{l}
{\it LL}\, := \,[(x_1)_{1}-{\it cx3}_{0} \, (u_1)_{0},\quad
(x_2)_{1}-{\it sx3}_{0} \, (u_1)_{0},\quad
(x_3)_{1}-(u_2)_{0},\quad
(y_1)_{0}-(x_1)_{0},\\[0.5em] \qquad
(y_2)_{0}-(x_2)_{0},\quad
{{\it cx3}_{0}}^{2}+{{\it sx3}_{0}}^{2}-1,\quad
{\it cx3}_{1}+{\it sx3}_{0} \, (x_3)_{1},\quad
{\it sx3}_{1}-{\it cx3}_{0} \, (x_3)_{1}]\\[0.2em]
\end{array}
\]}
\end{maplelatex}
\end{maplegroup}
\medskip
\noindent
We compute a Thomas decomposition with respect to $>$ of the given system of ordinary
differential equations.
\medskip

\begin{maplegroup}
\begin{mapleinput}
\mapleinline{active}{1d}{TD := DifferentialThomasDecomposition(LL, []);
}{}
\end{mapleinput}
\mapleresult
\begin{maplelatex}
\mapleinline{inert}{2d}{TD := [DifferentialSystem, DifferentialSystem, DifferentialSystem, DifferentialSystem, DifferentialSystem, DifferentialSystem, DifferentialSystem]}{
\[
\begin{array}{l}
{\it TD} := [{\it DifferentialSystem},{\it DifferentialSystem},{\it DifferentialSystem},\\[0.5em] \qquad
{\it DifferentialSystem},{\it DifferentialSystem},{\it DifferentialSystem},{\it DifferentialSystem}]\\[0.2em]
\end{array}
\]}
\end{maplelatex}
\end{maplegroup}
\medskip
\noindent
The first simple differential system is given as follows.
\medskip

\begin{maplegroup}
\begin{mapleinput}
\mapleinline{active}{1d}{Print(TD[1]);}{}
\end{mapleinput}
\mapleresult
\begin{maplelatex}
\mapleinline{inert}{2d}{[x1-y1 = 0, x2-y2 = 0, u1*cx3-y1[t] = 0, u1*sx3-y2[t] = 0, y1[t]^2*x3[t]+y2[t]^2*x3[t]-y1[t]*y2[t, t]+y2[t]*y1[t, t] = 0, u1^2-y1[t]^2-y2[t]^2 = 0, y1[t]^2*u2+y2[t]^2*u2-y1[t]*y2[t, t]+y2[t]*y1[t, t] = 0, y2[t] <> 0, y1[t] <> 0, y1[t]^2+y2[t]^2 <> 0, y2[t]*y1[t, t]-y1[t]*y2[t, t] <> 0]}{
\[
\begin{array}{l}
[\underline{x_1}-y_1=0,\quad
\underline{x_2}-y_2=0,\quad
u_1 \, \underline{{\it cx3}}-(y_1)_{t}=0,\quad
u_1 \, \underline{{\it sx3}}-(y_2)_{t}=0,\\[0.5em] \quad
(y_1)_{t}^2 \, \underline{(x_3)_{t}}+(y_2)_{t}^2 \, \underline{(x_3)_{t}}-(y_1)_{t} \, (y_2)_{t,t}+(y_2)_{t} \, (y_1)_{t,t}=0,\\[0.5em] \quad
\underline{u_1}^2-(y_1)_{t}^2-(y_2)_{t}^2=0,\quad
(y_1)_{t}^2 \, \underline{u_2}+(y_2)_{t}^2 \, \underline{u_2}-(y_1)_{t} \, (y_2)_{t,t}+(y_2)_{t} \, (y_1)_{t,t}=0,\\[0.5em] \quad
\underline{(y_2)_{t}}\neq 0,\quad
\underline{(y_1)_{t}}\neq 0,\quad
\underline{(y_1)_{t}}^2+{(y_2)_{t}}^2\neq 0,\quad
(y_2)_{t} \, \underline{(y_1)_{t,t}}-(y_1)_{t} \, (y_2)_{t,t}\neq 0]\\[0.2em]
\end{array}
\]}
\end{maplelatex}
\end{maplegroup}
\begin{maplegroup}
\begin{mapleinput}
\mapleinline{active}{1d}{collect(
}{}
\end{mapleinput}
\mapleresult
\begin{maplelatex}
\mapleinline{inert}{2d}{(y1[t]^2+y2[t]^2)*u2-y1[t]*y2[t, t]+y2[t]*y1[t, t] = 0}{
\[
\displaystyle  \left( (y_1)_{t}^2+(y_2)_{t}^2 \right) u_2-(y_1)_{t} \, (y_2)_{t,t}+(y_2)_{t} \, (y_1)_{t,t}=0
\]}
\end{maplelatex}
\end{maplegroup}
\medskip
\noindent
Thus, the equations with leader $u_1$ and $u_2$ in ${\it TD}[1]$ allow
to express $u_1$ and $u_2$ in terms of $y_1$ and $y_2$. (Up to solving
these equations for $u_1$ and $u_2$ explicitly, it is the same result
as in \cite[Example~5.12]{ConteMoogPerdon}.)
\medskip

The remaining six simple differential systems describe particular
configurations, which exhibit obstructions to invertibility.
\medskip

\begin{maplegroup}
\begin{mapleinput}
\mapleinline{active}{1d}{Print(TD[2]);}{}
\end{mapleinput}
\mapleresult
\begin{maplelatex}
\mapleinline{inert}{2d}{[x1-y1 = 0, x2-y2 = 0, u1*cx3-y1[t] = 0, u1*sx3-y2[t] = 0, x3[t] = 0, u1^2-y1[t]^2-y2[t]^2 = 0, u2 = 0, y2[t]*y1[t, t]-y1[t]*y2[t, t] = 0, y2[t] <> 0, y1[t] <> 0, y1[t]^2+y2[t]^2 <> 0]}{
\[
\begin{array}{l}
[\underline{x_1}-y_1=0,\quad
\underline{x_2}-y_2=0,\quad
u_1 \, \underline{{\it cx3}}-(y_1)_{t}=0,\quad
u_1 \, \underline{{\it sx3}}-(y_2)_{t}=0,\quad
(x_3)_{t}=0,\\[0.5em] \qquad
\underline{u_1}^2-(y_1)_{t}^2-(y_2)_{t}^2=0,\quad
\underline{u_2}=0,\quad
(y_2)_{t} \, \underline{(y_1)_{t,t}}-(y_1)_{t} \, (y_2)_{t,t}=0,\\[0.5em] \qquad
\underline{(y_2)_{t}}\neq 0,\quad
\underline{(y_1)_{t}}\neq 0,\quad
\underline{(y_1)_{t}}^2+{(y_2)_{t}}^2\neq 0]\\[0.2em]
\end{array}
\]}
\end{maplelatex}
\end{maplegroup}
\medskip
\noindent
The vanishing of the Wronskian determinant of $(y_1)_t$ and $(y_2)_t$
expresses that one of the velocities $\dot{x}_1$ and $\dot{x}_2$ is a
constant multiple of the other. Hence, no rotation is allowed, which forces
the input $u_2$ to be the zero function. Due to the inequations, the vector
$(\dot{x}_1, \dot{x}_2)$ is non-zero and not parallel to any of the $x_1$- or $x_2$-coordinate axes.
\medskip

\begin{maplegroup}
\begin{mapleinput}
\mapleinline{active}{1d}{Print(TD[3]);}{}
\end{mapleinput}
\mapleresult
\begin{maplelatex}
\mapleinline{inert}{2d}{[x1-y1 = 0, x2-y2 = 0, cx3+1 = 0, sx3 = 0, x3[t] = 0, u1+y1[t] = 0, u2 = 0, y2[t] = 0, y1[t] <> 0]}{
\[
\begin{array}{l}
[\underline{x_1}-y_1=0,\quad
\underline{x_2}-y_2=0,\quad
\underline{{\it cx3}}+1=0,\quad
\underline{{\it sx3}}=0,\quad
\underline{(x_3)_{t}}=0,\\[0.5em] \qquad
\underline{u_1}+(y_1)_{t}=0,\quad
\underline{u_2}=0,\quad
\underline{(y_2)_{t}}=0,\quad
\underline{(y_1)_{t}}\neq 0]\\[0.2em]
\end{array}
\]}
\end{maplelatex}
\end{maplegroup}
\begin{maplegroup}
\begin{mapleinput}
\mapleinline{active}{1d}{Print(TD[4]);}{}
\end{mapleinput}
\mapleresult
\begin{maplelatex}
\mapleinline{inert}{2d}{[x1-y1 = 0, x2-y2 = 0, cx3-1 = 0, sx3 = 0, x3[t] = 0, u1-y1[t] = 0, u2 = 0, y2[t] = 0, y1[t] <> 0]}{
\[
\begin{array}{l}
[\underline{x_1}-y_1=0,\quad
\underline{x_2}-y_2=0,\quad
\underline{{\it cx3}}-1=0,\quad
\underline{{\it sx3}}=0,\quad
\underline{(x_3)_{t}}=0,\\[0.5em] \qquad
\underline{u_1}-(y_1)_{t}=0,\quad
\underline{u_2}=0,\quad
\underline{(y_2)_{t}}=0,\quad
\underline{(y_1)_{t}}\neq 0]\\[0.2em]
\end{array}
\]}
\end{maplelatex}
\end{maplegroup}
\medskip
\noindent
The previous two simple systems describe cases in which only
movement in any of the two directions defined by the $x_1$-coordinate axis
is allowed and no rotation.
\medskip

\begin{maplegroup}
\begin{mapleinput}
\mapleinline{active}{1d}{Print(TD[5]);}{}
\end{mapleinput}
\mapleresult
\begin{maplelatex}
\mapleinline{inert}{2d}{[x1-y1 = 0, x2-y2 = 0, cx3^2+sx3^2-1 = 0, sx3[t]-u2*cx3 = 0, x3[t]-u2 = 0, u1 = 0, y1[t] = 0, y2[t] = 0, sx3+1 <> 0, sx3-1 <> 0]}{
\[
\begin{array}{l}
[\underline{x_1}-y_1=0,\quad
\underline{x_2}-y_2=0,\quad
\underline{{{\it cx3}}}^{2}+{{\it sx3}}^{2}-1=0,\quad
\underline{{\it sx3}_{t}}-u_2\,{\it cx3}=0,\\[0.5em] \quad
\underline{(x_3)_{t}}-u_2=0,\quad
\underline{u_1}=0,\quad
\underline{(y_1)_{t}}=0,\quad
\underline{(y_2)_{t}}=0,\quad
\underline{{\it sx3}}+1\neq 0,\quad
\underline{{\it sx3}}-1\neq 0]\\[0.2em]
\end{array}
\]}
\end{maplelatex}
\end{maplegroup}
\medskip
\noindent
The fifth simple system describes configurations which only allow rotation,
and the input $u_1$ is forced to be the zero function.
(Similarly to Example~\ref{mlhdr:ex:thomascircle},
the inequation ${\it sx3}^2 - 1 \neq 0$ is introduced here to
ensure that ${\it cx3}^2 + {\it sx3}^2 - 1$ has no multiple roots
as polynomial in ${\it cx3}$. It is included in the simple system
in factorized form.)
\medskip

\begin{maplegroup}
\begin{mapleinput}
\mapleinline{active}{1d}{Print(TD[6]);}{}
\end{mapleinput}
\mapleresult
\begin{maplelatex}
\mapleinline{inert}{2d}{[x1-y1 = 0, x2-y2 = 0, cx3 = 0, sx3+1 = 0, x3[t] = 0, u1+y2[t] = 0, u2 = 0, y1[t] = 0]}{
\[
\begin{array}{l}
[\underline{x_1}-y_1=0,\quad
\underline{x_2}-y_2=0,\quad
\underline{{\it cx3}}=0,\quad
\underline{{\it sx3}}+1=0,\quad
\underline{(x_3)_{t}}=0,\\[0.5em] \qquad
\underline{u_1}+(y_2)_t=0,\quad
\underline{u_2}=0,\quad
\underline{(y_1)_{t}}=0]\\[0.2em]
\end{array}
\]}
\end{maplelatex}
\end{maplegroup}
\begin{maplegroup}
\begin{mapleinput}
\mapleinline{active}{1d}{Print(TD[7]);}{}
\end{mapleinput}
\mapleresult
\begin{maplelatex}
\mapleinline{inert}{2d}{[x1-y1 = 0, x2-y2 = 0, cx3 = 0, sx3-1 = 0, x3[t] = 0, u1-y2[t] = 0, u2 = 0, y1[t] = 0]}{
\[
\begin{array}{l}
[\underline{x_1}-y_1=0,\quad
\underline{x_2}-y_2=0,\quad
\underline{{\it cx3}}=0,\quad
\underline{{\it sx3}}-1=0,\quad
\underline{(x_3)_{t}}=0,\\[0.5em] \qquad
\underline{u_1}-(y_2)_{t}=0,\quad
\underline{u_2}=0,\quad
\underline{(y_1)_{t}}=0]\\[0.2em]
\end{array}
\]}
\end{maplelatex}
\end{maplegroup}
\medskip
\noindent
The last two simple systems cover the cases of movement in any of
the two directions defined by the $x_2$-coordinate axis and no rotation.
\end{example}

\medskip

Next we consider the detection of flat outputs.

\begin{example}
A model of a 2-D crane is given by the following system of
ordinary differential equations (cf.\ \cite[Sect.~4.1]{FliessLevineMartinRouchon} and
the references therein), where $x(t)$ and $z(t)$ are the coordinates of the
load of mass $m$, $\theta(t)$ is the angle between the rope and the
$z$-axis, $d(t)$ the trolley position, $T(t)$ the tension of the rope,
$R(t)$ the rope length, and $g$ the gravitational constant.
\[
\left\{
\begin{array}{rcl}
m \, \ddot{x} & = & -T \, \sin \theta,\\[0.5em]
m \, \ddot{z} & = & -T \, \cos \theta + m \, g,\\[0.5em]
x & = & R \, \sin \theta + d,\\[0.5em]
z & = & R \, \cos \theta.
\end{array}
\right.
\]
The task is to decide whether $\{ x, z \}$ is a flat output of the system.

\medskip

Similarly to the previous example, we represent $\cos \theta$ and $\sin \theta$
by differential indeterminates $c$ and $s$ and add the generating relation
$c^2 + s^2 = 1$ to the system.
In this example the given equations depend on $\theta$ only
through $\cos \theta$ and $\sin \theta$. Therefore, we do not include $\theta$
as a differential indeterminate and do not need to
add the relations $c_t = -s \, \theta_t$ and $s_t = c \, \theta_t$ to the system.
(Note that, if $I$ is the differential ideal of $\mathbb{Q}\{ \theta, c, s \}$ with
derivation $\partial_t$ which is generated by $c^2 + s^2 - 1$
and $c_t + s \, \theta_t$ and $s_t - c \, \theta_t$, then $I \cap \mathbb{Q}\{ c, s \}$
is the differential ideal which is generated by $c^2 + s^2 - 1$.)
\bigskip

\begin{maplegroup}
\begin{mapleinput}
\mapleinline{active}{1d}{with(DifferentialThomas):}{}
\end{mapleinput}
\end{maplegroup}
\begin{maplegroup}
\begin{mapleinput}
\mapleinline{active}{1d}{ivar := [t]:}{}
\end{mapleinput}
\end{maplegroup}
\begin{maplegroup}
\begin{mapleinput}
\mapleinline{active}{1d}{dvar := [T,c,s,d,R,x,z]:}{}
\end{mapleinput}
\end{maplegroup}
\medskip
\noindent
We set up the block ranking $>$ which satisfies $\{ T, c, s, d, R \} \gg \{ x, z \}$
as well as $T > c > s > d > R$ and $x > z$.
\medskip

\begin{maplegroup}
\begin{mapleinput}
\mapleinline{active}{1d}{ComputeRanking(ivar, [[T,c,s,d,R],[x,z]]):}{}
\end{mapleinput}
\end{maplegroup}
\medskip
\noindent
We compute a Thomas decomposition with respect to $>$.
(As is customary in Maple, the symbols $m$ and $g$ are treated here as
algebraically independent over $\mathbb{Q}$. More precisely, the ground
field for the following computation is the differential field $\mathbb{Q}(m, g)$
with trivial derivation.)
\medskip

\begin{maplegroup}
\begin{mapleinput}
\mapleinline{active}{1d}{TD := DifferentialThomasDecomposition(
[m*x[2]+T[0]*s[0], m*z[2]+T[0]*c[0]-m*g,
x[0]-R[0]*s[0]-d[0], z[0]-R[0]*c[0],
c[0]\symbol{94}2+s[0]\symbol{94}2-1], []);
}{}
\end{mapleinput}
\mapleresult
\begin{maplelatex}
\mapleinline{inert}{2d}{TD := [DifferentialSystem, DifferentialSystem, DifferentialSystem, DifferentialSystem, DifferentialSystem, DifferentialSystem, DifferentialSystem]}{\[
\begin{array}{l}
{\it TD} := [{\it DifferentialSystem},{\it DifferentialSystem},{\it DifferentialSystem},\\[0.5em] \qquad
{\it DifferentialSystem},{\it DifferentialSystem},{\it DifferentialSystem},\\[0.5em] \qquad
{\it DifferentialSystem}]\\[0.2em]
\end{array}
\]}
\end{maplelatex}
\end{maplegroup}
\medskip
\noindent
The second simple differential system is given as follows.
\medskip

\begin{maplegroup}
\begin{mapleinput}
\mapleinline{active}{1d}{Print(TD[2]);
}{}
\end{mapleinput}
\mapleresult
\begin{maplelatex}
\mapleinline{inert}{2d}{[z*T+m*z[t, t]*R-m*g*R = 0, R*c-z = 0, z[t, t]*R*s-g*R*s-z*x[t, t] = 0, z[t, t]*d-g*d+z*x[t, t]-x*z[t, t]+g*x = 0, z[t, t]^2*R^2-2*g*z[t, t]*R^2+g^2*R^2-z^2*x[t, t]^2-z^2*z[t, t]^2+2*g*z^2*z[t, t]-g^2*z^2 = 0, z <> 0, z[t, t]-g <> 0, x[t, t] <> 0, x[t, t]^2+z[t, t]^2-2*g*z[t, t]+g^2 <> 0]}{
\[
\begin{array}{l}
[z \, \underline{T}+m \, z_{{t,t}} \, R-m \, g \, R=0,\quad
R \, \underline{c}-z=0,\quad
z_{t,t} \, R \, \underline{s}-g \, R \, \underline{s}-z \, x_{t,t}=0,\\[0.5em] \qquad
z_{t,t} \, \underline{d}-g \, \underline{d}+z \, x_{t,t}-x \, z_{t,t}+g \, x=0,\\[0.5em] \qquad
{z_{t,t}}^2 \, \underline{R}^2-2 \, g \, z_{t,t} \, \underline{R}^2+g^2 \, \underline{R}^2-z^2 \, x_{t,t}^2-z^2 \, z_{t,t}^2+2 \, g \, z^2 \, z_{t,t}-g^2 \, z^2=0,\\[0.5em] \qquad
\underline{z}\neq 0,\quad
\underline{z_{t,t}}-g\neq 0,\quad
\underline{x_{t,t}}\neq 0,\quad
\underline{{x_{t,t}}}^2+{z_{t,t}}^{2}-2 \, g \, z_{t,t}+g^2\neq 0]\\[0.2em]
\end{array}
\]}
\end{maplelatex}
\end{maplegroup}
\begin{maplegroup}
\begin{mapleinput}
\mapleinline{active}{1d}{collect(
}{}
\end{mapleinput}
\mapleresult
\begin{maplelatex}
\mapleinline{inert}{2d}{(z[t, t]-g)^2*R^2-z^2*(z[t, t]^2-2*g*z[t, t]+g^2+x[t, t]^2) = 0}{\[\displaystyle  \left( z_{{t,t}}-g \right) ^{2}{R}^{2}-{z}^{2} \left( {x_{{t,t}}}^{2}+{z_{{t,t}}}^{2}-2 \, g \, z_{{t,t}}+{g}^{2} \right) =0\]}
\end{maplelatex}
\end{maplegroup}
\medskip
\noindent
We observe that this simple system $S$ contains no equation involving
derivatives of $x$ and $z$ only. Moreover, the equations in $S$ show
that $T$, $c$, $s$, $d$, $R$ are observable with respect to $\{ x, z \}$.
Hence, $\{ x, z \}$ is a flat output of $S$.
\medskip

The remaining six simple differential systems describe particular configurations
for which $\{ x, z \}$ is not a flat output. In fact, the movement of the load
is restricted by some constraint in these cases (e.g., $x_{t,t} = 0$ or $z = 0$,
one reason being, e.g., that vanishing rope tension implies constant acceleration
of the load, another being a constant rope length of zero allowing no vertical
movement of the load).
We do not consider the system to be controllable under these conditions.
\medskip

\begin{maplegroup}
\begin{mapleinput}
\mapleinline{active}{1d}{Print(TD[1]);
}{}
\end{mapleinput}
\mapleresult
\begin{maplelatex}
\mapleinline{inert}{2d}{[T = 0, R*c-z = 0, R*s+d-x = 0, d^2-2*x*d+x^2-R^2+z^2 = 0, x[t, t] = 0, z[t, t]-g = 0, z <> 0, R <> 0, R+z <> 0, R-z <> 0]}{
\[
\begin{array}{l}
[\underline{T}=0,\quad
R \, \underline{c}-z=0,\quad
R \, \underline{s}+d-x=0,\quad
\underline{d}^2-2 \, x \, \underline{d}+x^2-R^2+z^2=0,\\[0.5em] \qquad
\underline{x_{{t,t}}}=0,\quad
\underline{z_{{t,t}}}-g=0,\quad
\underline{z}\neq 0,\quad
\underline{R}\neq 0,\quad
\underline{R}+z\neq 0,\quad
\underline{R}-z\neq 0]\\[0.2em]
\end{array}
\]}
\end{maplelatex}
\end{maplegroup}
\begin{maplegroup}
\begin{mapleinput}
\mapleinline{active}{1d}{Print(TD[3]);
}{}
\end{mapleinput}
\mapleresult
\begin{maplelatex}
\mapleinline{inert}{2d}{[T-m*z[t, t]+m*g = 0, c+1 = 0, s = 0, d-x = 0, R+z = 0, x[t, t] = 0, z <> 0]}{
\[
\begin{array}{l}
[\underline{T}-m \, z_{t,t}+m \, g=0,\quad
\underline{c}+1=0,\quad
\underline{s}=0,\quad
\underline{d}-x=0,\quad
\underline{R}+z=0,\\[0.5em] \qquad
\underline{x_{{t,t}}}=0,\quad
\underline{z}\neq 0]\\[0.2em]
\end{array}
\]}
\end{maplelatex}
\end{maplegroup}
\begin{maplegroup}
\begin{mapleinput}
\mapleinline{active}{1d}{Print(TD[4]);
}{}
\end{mapleinput}
\mapleresult
\begin{maplelatex}
\mapleinline{inert}{2d}{[T+m*z[t, t]-m*g = 0, c-1 = 0, s = 0, d-x = 0, R-z = 0, x[t, t] = 0, z <> 0]}{
\[
\begin{array}{l}
[\underline{T}+m \, z_{t,t}-m \, g=0,\quad
\underline{c}-1=0,\quad
\underline{s}=0,\quad
\underline{d}-x=0,\quad
\underline{R}-z=0,\\[0.5em] \qquad
\underline{x_{t,t}}=0,\quad
\underline{z}\neq 0]\\[0.2em]
\end{array}
\]}
\end{maplelatex}
\end{maplegroup}
\begin{maplegroup}
\begin{mapleinput}
\mapleinline{active}{1d}{Print(TD[5]);
}{}
\end{mapleinput}
\mapleresult
\begin{maplelatex}
\mapleinline{inert}{2d}{[s*T+m*x[t, t] = 0, x[t, t]*c+g*s = 0, g^2*s^2+x[t, t]^2*s^2-x[t, t]^2 = 0, d-x = 0, R = 0, z = 0, x[t, t] <> 0, x[t, t]^2+g^2 <> 0]}{
\[
\begin{array}{l}
[s \, \underline{T}+m \, x_{t,t}=0,\quad
x_{t,t} \, \underline{c}+g \, s=0,\quad
g^2 \, \underline{s}^2+x_{t,t}^2 \, \underline{s}^2-x_{t,t}^2=0,\quad
\underline{d}-x=0,\quad
\underline{R}=0,\\[0.5em] \qquad
\underline{z}=0,\quad
\underline{x_{t,t}}\neq 0,\quad
\underline{{x_{t,t}}}^2+g^2\neq 0]\\[0.2em]
\end{array}
\]}
\end{maplelatex}
\end{maplegroup}
\begin{maplegroup}
\begin{mapleinput}
\mapleinline{active}{1d}{Print(TD[6]);
}{}
\end{mapleinput}
\mapleresult
\begin{maplelatex}
\mapleinline{inert}{2d}{[T+m*g = 0, c+1 = 0, s = 0, d-x = 0, R = 0, x[t, t] = 0, z = 0]}{
\[
[\underline{T}+m \, g=0,\quad
\underline{c}+1=0,\quad
\underline{s}=0,\quad
\underline{d}-x=0,\quad
\underline{R}=0,\quad
\underline{x_{t,t}}=0,\quad
\underline{z}=0]
\]}
\end{maplelatex}
\end{maplegroup}
\begin{maplegroup}
\begin{mapleinput}
\mapleinline{active}{1d}{Print(TD[7]);
}{}
\end{mapleinput}
\mapleresult
\begin{maplelatex}
\mapleinline{inert}{2d}{[\underline{T}-m*g = 0, c-1 = 0, s = 0, d-x = 0, R = 0, x[t, t] = 0, z = 0]}{
\[
\displaystyle [\underline{T}-m \, g=0,\quad
\underline{c}-1=0,\quad
\underline{s}=0,\quad
\underline{d}-x=0,\quad
\underline{R}=0,\quad
\underline{x_{t,t}}=0,\quad
\underline{z}=0]
\]}
\end{maplelatex}
\end{maplegroup}
\end{example}

\medskip

We give two examples which demonstrate how the Thomas decomposition
technique can be used to study the dependence of structural properties
of a nonlinear control system on parameters.

\begin{example}
A model of a continuous stirred-tank reactor (cf.\ \cite[Example~1.2]{KwakernaakSivan})
is given by the differential system 
\[
\left\{
\begin{array}{rcl}
\dot{V}(t) & = & F_1(t) + F_2(t) - k \, \sqrt{V(t)},\\[0.5em]
\dot{\overline{c(t) \, V(t)}} & = & c_1 \, F_1(t) + c_2 \, F_2(t) - c(t) \, k \, \sqrt{V(t)}.
\end{array}
\right.
\]
A dissolved material has concentration $c(t)$ in the tank and it is fed
through two inputs with constant concentrations $c_1$ and $c_2$ and flow
rates $F_1(t)$ and $F_2(t)$, respectively. There exists an outward flow with a
flow rate proportional to the square root of the volume $V(t)$ of liquid in
the tank. Moreover, $k$ is an experimental constant.

\medskip

In order to eliminate the square root of the volume in the given equations,
we represent $\sqrt{V(t)}$ as a differential indeterminate $sV$ and substitute
other occurrences of $V(t)$ by $sV^2$. We investigate the dependence of the
behavior on parameter configurations by considering $c_1$ and $c_2$ as
differential indeterminates as well and adding the conditions $\dot{c}_1 = 0$
and $\dot{c}_2 = 0$.

\bigskip

\begin{maplegroup}
\begin{mapleinput}
\mapleinline{active}{1d}{with(DifferentialThomas):
}{}
\end{mapleinput}
\end{maplegroup}
\begin{maplegroup}
\begin{mapleinput}
\mapleinline{active}{1d}{ivar := [t]:
}{}
\end{mapleinput}
\end{maplegroup}
\begin{maplegroup}
\begin{mapleinput}
\mapleinline{active}{1d}{dvar := [F1,F2,sV,c,c1,c2]:
}{}
\end{mapleinput}
\end{maplegroup}
\medskip
\noindent
We define $R = \mathbb{Q}\{ F_1, F_2, sV, c, c_1, c_2 \}$
and choose the block ranking $>$ on $R$ with blocks $\{ F_1, F_2 \}$,
$\{ sV, c \}$, $\{ c_1, c_2 \}$, i.e., satisfying
$\{ F_2, F_2 \} \gg \{ sV, c \} \gg \{ c_1, c_2 \}$
and $F_1 > F_2$ and $sV > c$ and $c_1 > c_2$.
\medskip

\begin{maplegroup}
\begin{mapleinput}
\mapleinline{active}{1d}{ComputeRanking(ivar, [[F1,F2],[sV,c],[c1,c2]]):
}{}
\end{mapleinput}
\end{maplegroup}
\begin{maplegroup}
\begin{mapleinput}
\mapleinline{active}{1d}{L := [2*sV[t]*sV-F1-F2+k*sV,
c[t]*sV\symbol{94}2-c2*F2+c*k*sV-c1*F1+2*c*sV[t]*sV,
c1[t], c2[t]]:
}{}
\end{mapleinput}
\end{maplegroup}
\begin{maplegroup}
\begin{mapleinput}
\mapleinline{active}{1d}{LL := Diff2JetList(Ind2Diff(L, ivar, dvar));
}{}
\end{mapleinput}
\mapleresult
\begin{maplelatex}
\mapleinline{inert}{2d}{LL := [2*sV[1]*sV[0]-F1[0]-F2[0]+k*sV[0], c[1]*sV[0]^2-c2[0]*F2[0]+c[0]*k*sV[0]-c1[0]*F1[0]+2*c[0]*sV[1]*sV[0], c1[1], c2[1]]}{
\[
\begin{array}{l}
{\it LL}\, := \,[2 \, {\it sV}_{1} \, {\it sV}_{0}-(F_1)_{0}-(F_2)_{0}+k \, {\it sV}_{0},\\[0.5em] \qquad
c_{1} \, {{\it sV}_{0}}^2-(c_2)_{0} \, (F_2)_{0}+c_{0} \, k \, {\it sV}_{0}-(c_1)_{0} \, (F_1)_{0}+2 \, c_{0} \, {\it sV}_{1} \, {\it sV}_{0},\quad
(c_1)_{1},\quad
(c_2)_{1}]\\[0.2em]
\end{array}
\]}
\end{maplelatex}
\end{maplegroup}
\medskip
\noindent
We compute a Thomas decomposition with respect to $>$ of the given system of ordinary
differential equations, to which we add the inequations $\sqrt{V} \neq 0$,
$c_1 \neq 0$, $c_2 \neq 0$ to exclude trivial cases.
\medskip

\begin{maplegroup}
\begin{mapleinput}
\mapleinline{active}{1d}{TD := DifferentialThomasDecomposition(LL,
[sV[0],c1[0],c2[0]]);
}{}
\end{mapleinput}
\mapleresult
\begin{maplelatex}
\mapleinline{inert}{2d}{TD := [DifferentialSystem, DifferentialSystem, DifferentialSystem]}{\[\displaystyle {\it TD}\, := \,[{\it DifferentialSystem},{\it DifferentialSystem},{\it DifferentialSystem}]\]}
\end{maplelatex}
\end{maplegroup}
\medskip
\noindent
The first simple differential system is given as follows.
\medskip

\begin{maplegroup}
\begin{mapleinput}
\mapleinline{active}{1d}{Print(TD[1]);
}{}
\end{mapleinput}
\mapleresult
\begin{maplelatex}
\mapleinline{inert}{2d}{[c2*F1-c1*F1+2*c*sV*sV[t]-2*c2*sV*sV[t]+c[t]*sV^2+c*k*sV-c2*k*sV = 0, c1*F2-c2*F2+2*c*sV*sV[t]-2*c1*sV*sV[t]+c[t]*sV^2+c*k*sV-c1*k*sV = 0, c1[t] = 0, c2[t] = 0, c2 <> 0, c1 <> 0, c1-c2 <> 0, sV <> 0]}{
\[
\begin{array}{l}
[c_2 \, \underline{F_1}-c_1 \, \underline{F_1}+2 \, c \, {\it sV}{\it sV}_{t}-2 \, c_2 \, {\it sV} \, {\it sV}_{t}+c_{t} \, {{\it sV}}^2+c \, k \, {\it sV}-c_2 \, k \, {\it sV}=0,\\[0.5em] \qquad
c_1 \, \underline{F_2}-c_2 \, \underline{F_2}+2 \, c \, {\it sV} \, {\it sV}_{t}-2 \, c_1 \, {\it sV} \, {\it sV}_{t}+c_{t} \, {{\it sV}}^2+c \, k \, {\it sV}-c_1 \, k \, {\it sV}=0,\\[0.5em] \qquad
\underline{(c_1)_{t}}=0,\quad
\underline{(c_2)_{t}}=0,\quad
\underline{c_2}\neq 0,\quad
\underline{c_1}\neq 0,\quad
\underline{c_1}-c_2\neq 0,\quad
\underline{{\it sV}}\neq 0]\\[0.2em]
\end{array}
\]}
\end{maplelatex}
\end{maplegroup}
\begin{maplegroup}
\begin{mapleinput}
\mapleinline{active}{1d}{collect(
}{}
\end{mapleinput}
\mapleresult
\begin{maplelatex}
\mapleinline{inert}{2d}{(c2-c1)*F1+2*c*sV*sV[t]-2*c2*sV*sV[t]+c[t]*sV^2+c*k*sV-c2*k*sV = 0}{
\[
\displaystyle \left( c_2-c_1 \right) F_1+2 \, c \, {\it sV} \, {\it sV}_{t}-2 \, c_2 \, {\it sV} \, {\it sV}_{t}+c_{t} \, {\it sV}^2+c \, k \, {\it sV}-c_2 \, k \, {\it sV}=0
\]}
\end{maplelatex}
\end{maplegroup}
\begin{maplegroup}
\begin{mapleinput}
\mapleinline{active}{1d}{collect(
}{}
\end{mapleinput}
\mapleresult
\begin{maplelatex}
\mapleinline{inert}{2d}{(c1-c2)*F2+2*c*sV*sV[t]-2*c1*sV*sV[t]+c[t]*sV^2+c*k*sV-c1*k*sV = 0}{
\[
\displaystyle \left( c_1-c_2 \right) F_2+2 \, c \, {\it sV} \, {\it sV}_{t}-2 \, c_1 \, {\it sV} \, {\it sV}_{t}+c_{t} \, {\it sV}^2+c \, k \, {\it sV}-{\it c1} \, k \, {\it sV}=0
\]}
\end{maplelatex}
\end{maplegroup}
\medskip
\noindent
The first two equations in the first simple system $S$ show that
$F_1$ and $F_2$ are observable with respect to $\{ c, sV \}$. (Although $c_1$
and $c_2$ are represented by differential indeterminates here, we consider
these still as parameters.)
Let $E$ be the differential ideal of $R$ generated by $S^{=}$ and
$q$ the product of the initials (and separants) of all elements of $S^{=}$.
Due to the choice of the block ranking, we conclude that we have
$(E : q^{\infty}) \cap \mathbb{Q}\{ sV, c \} = \{ 0 \}$ (cf.\
Proposition~\ref{mlhdr:prop:elimsimple}).
Hence, $\{ c, sV \}$ is a flat output of $S$.
\medskip

The remaining two simple systems describe configurations of the system
in which the two concentrations $c_1$ and $c_2$ are equal. Since both
input feeds are identical and constant, this condition precludes control of the
concentration in the tank. These particular systems do not admit $\{ c, sV \}$
as a flat output. In fact, by inspecting the equations of these systems,
we observe that we have $(E : q^{\infty}) \cap \mathbb{Q}\{ sV, c \} \neq \{ 0 \}$.

\medskip

\begin{maplegroup}
\begin{mapleinput}
\mapleinline{active}{1d}{Print(TD[2]);
}{}
\end{mapleinput}
\mapleresult
\begin{maplelatex}
\mapleinline{inert}{2d}{[c*F1-c2*F1+c*F2-c2*F2+c[t]*sV^2 = 0, 2*c*sV[t]-2*c2*sV[t]+c[t]*sV+c*k-c2*k = 0, c1-c2 = 0, c2[t] = 0, c2 <> 0, c-c2 <> 0, sV <> 0]}{
\[
\begin{array}{l}
[c \, \underline{F_1}-c_2\,\underline{F_1}+c \, F_2-c_2 \, F_2+c_{t} \, {\it sV}^2=0,\\[0.5em] \qquad
2 \, c \, \underline{{\it sV}_{{t}}}-2 \, c_2 \, \underline{{\it sV}_{t}}+c_{t} \, {\it sV}+c \, k-c_2 \, k=0,\quad
\underline{c_1}-c_2=0,\quad
\underline{(c_2)_{t}}=0,\\[0.5em] \qquad
\underline{c_2}\neq 0,\quad
\underline{c}-c_2\neq 0,\quad
\underline{{\it sV}}\neq 0]\\[0.2em]
\end{array}
\]}
\end{maplelatex}
\end{maplegroup}
\begin{maplegroup}
\begin{mapleinput}
\mapleinline{active}{1d}{Print(TD[3]);
}{}
\end{mapleinput}
\mapleresult
\begin{maplelatex}
\mapleinline{inert}{2d}{[F1+F2-2*sV*sV[t]-k*sV = 0, c-c2 = 0, c1-c2 = 0, c2[t] = 0, c2 <> 0, sV <> 0]}{
\[
\begin{array}{l}
[\underline{F_1}+F_2-2 \, {\it sV} \, {\it sV}_{t}-k \, {\it sV}=0,\quad
\underline{c}-c_2=0,\quad
\underline{c_1}-c_2=0,\quad
\underline{(c_2)_{t}}=0,\\[0.5em] \qquad
\underline{c_2}\neq 0,\quad
\underline{{\it sV}}\neq 0]\\[0.2em]
\end{array}
\]}
\end{maplelatex}
\end{maplegroup}
\end{example}

\begin{example}
Let us consider the following system of linear partial differential
equations for functions $\xi_1$, $\xi_2$, $\xi_3$ of
$\mathbf{x} = (x_1, x_2, x_3)$
involving a parametric function $a(x_2)$
\[
\left\{
\begin{array}{rcl}
\displaystyle
-a(x_2) {\frac{\partial \xi_1(\mathbf{x})}{\partial x_1}}+
{\frac{\partial \xi_3(\mathbf{x})}{\partial x_1}} - \left( {\frac{\partial}{\partial x_2}}a(x_2)  \right) \xi_2(\mathbf{x})
+ \frac{1}{2}\,a(x_2) \left(\nabla\cdot\xi(\mathbf{x})\right) & = & 0,\\[1em]
\displaystyle
-a(x_2) {\frac{\partial \xi_1(\mathbf{x})}{\partial x_2}}+
{\frac{\partial \xi_3(\mathbf{x})}{\partial x_2}} & = & 0,\\[1em]
\displaystyle
-a(x_2) {\frac{\partial \xi_1(\mathbf{x})}{\partial x_3}}+
{\frac{\partial \xi_3(\mathbf{x})}{\partial x_3}}-
\frac{1}{2} \left(\nabla\cdot\xi(\mathbf{x})\right) & = & 0,
\end{array}
\right.
\]
which describe infinitesimal transformations associated to a certain
Pfaffian system \cite[Example~4]{PommaretQuadrat1997}. In order to
study the influence of the parametric function $a$ on the system
using the package {\tt DifferentialThomas}, $a$ is included in the
list of dependent variables and its dependence on merely $x_2$ is
taken into account by adding the following two equations to the system:
\[
{\frac{\partial}{\partial x_1}}a(x_1,x_2,x_3)=0,\qquad
{\frac{\partial}{\partial x_3}}a(x_1,x_2,x_3)=0.
\]
Let $R$ be the differential polynomial ring $\mathbb{Q}\{ \xi_1, \xi_2, \xi_3, a \}$,
endowed with the partial differential operators $\partial_1$, $\partial_2$,
$\partial_3$ with respect to $x_1$, $x_2$, $x_3$.
\bigskip

\begin{maplegroup}
\begin{mapleinput}
\mapleinline{active}{1d}{with(DifferentialThomas):
}{}
\end{mapleinput}
\end{maplegroup}
\begin{maplegroup}
\begin{mapleinput}
\mapleinline{active}{1d}{ivar := [x1,x2,x3]:}{}
\end{mapleinput}
\end{maplegroup}
\begin{maplegroup}
\begin{mapleinput}
\mapleinline{active}{1d}{dvar := [xi1,xi2,xi3,a]:}{}
\end{mapleinput}
\end{maplegroup}
\medskip
\noindent
We choose a block ranking $>$ on $R$ with blocks $\{ \xi_1, \xi_2, \xi_3 \}$,
$\{ a \}$.
\medskip

\begin{maplegroup}
\begin{mapleinput}
\mapleinline{active}{1d}{ComputeRanking(ivar, [[xi1,xi2,xi3],[a]]):}{}
\end{mapleinput}
\end{maplegroup}
\begin{maplegroup}
\begin{mapleinput}
\mapleinline{active}{1d}{L := [-a*xi1[x1]+xi3[x1]-a[x2]*xi2
+(1/2)*a*(xi1[x1]+xi2[x2]+xi3[x3]),
-a*xi1[x2]+xi3[x2], -a*xi1[x3]+xi3[x3]
-(1/2)*(xi1[x1]+xi2[x2]+xi3[x3]), a[x1], a[x3]]:
}{}
\end{mapleinput}
\end{maplegroup}
\begin{maplegroup}
\begin{mapleinput}
\mapleinline{active}{1d}{LL := Diff2JetList(Ind2Diff(L, ivar, dvar));
}{}
\end{mapleinput}
\mapleresult
\begin{maplelatex}
\mapleinline{inert}{2d}{LL := [-a[0, 0, 0]*xi1[1, 0, 0]+xi3[1, 0, 0]+(1/2)*a[0, 0, 0]*(xi1[1, 0, 0]+xi2[0, 1, 0]+xi3[0, 0, 1])-a[0, 1, 0]*xi2[0, 0, 0], -a[0, 0, 0]*xi1[0, 1, 0]+xi3[0, 1, 0], -a[0, 0, 0]*xi1[0, 0, 1]+(1/2)*xi3[0, 0, 1]-(1/2)*xi1[1, 0, 0]-(1/2)*xi2[0, 1, 0], a[1, 0, 0], a[0, 0, 1]]}{
\[
\begin{array}{l}
{\it LL}\, := \,[-a_{{0,0,0}} \, (\xi_1)_{{1,0,0}}+(\xi_3)_{{1,0,0}}+\frac{1}{2}\,a_{{0,0,0}} \left( (\xi_1)_{{1,0,0}}+(\xi_2)_{{0,1,0}}+(\xi_3)_{{0,0,1}} \right)\\[0.5em] \qquad
-a_{{0,1,0}} \, (\xi_2)_{{0,0,0}},\quad
-a_{{0,0,0}} \, (\xi_1)_{{0,1,0}}+(\xi_3)_{{0,1,0}},\\[0.5em] \qquad
-a_{{0,0,0}} \, (\xi_1)_{{0,0,1}}+\frac{1}{2}\,(\xi_3)_{{0,0,1}}-\frac{1}{2}\,(\xi_1)_{{1,0,0}}-\frac{1}{2}\,(\xi_2)_{{0,1,0}},\quad
a_{{1,0,0}},\quad
a_{{0,0,1}}]\\[0.2em]
\end{array}
\]}
\end{maplelatex}
\end{maplegroup}
\medskip
\noindent
We compute a Thomas decomposition with respect to $>$ of the given system
of partial differential equations.
\medskip

\begin{maplegroup}
\begin{mapleinput}
\mapleinline{active}{1d}{TD := DifferentialThomasDecomposition(LL, []);
}{}
\end{mapleinput}
\mapleresult
\begin{maplelatex}
\mapleinline{inert}{2d}{TD := [DifferentialSystem, DifferentialSystem, DifferentialSystem]}{
\[
\displaystyle {\it TD}\, := \,[{\it DifferentialSystem},{\it DifferentialSystem},{\it DifferentialSystem}]
\]}
\end{maplelatex}
\end{maplegroup}
\medskip
\noindent
The resulting three simple differential systems are given as follows.
\medskip

\begin{maplegroup}
\begin{mapleinput}
\mapleinline{active}{1d}{Print(TD[1]);
}{}
\end{mapleinput}
\mapleresult
\begin{maplelatex}
\mapleinline{inert}{2d}{[a*xi1[x2]-xi3[x2] = 0, a^2*xi1[x3]+xi3[x1] = 0, xi2 = 0, a*xi1[x1]-2*xi3[x1]-a*xi3[x3] = 0, a[x1] = 0, a[x3] = 0, a <> 0]}{
\[
\begin{array}{l}
[a \, \underline{(\xi_1)_{x_2}}-(\xi_3)_{x_2}=0,\quad
{a}^{2} \, \underline{(\xi_1)_{x_3}}+(\xi_3)_{x_1}=0,\quad
\underline{\xi_2}=0,\\[0.5em] \qquad
a \, \underline{(\xi_1)_{x_1}}-2\,(\xi_3)_{x_1}-a \, (\xi_3)_{x_3}=0,\quad
\underline{a_{x_1}}=0,\quad
\underline{a_{x_3}}=0,\quad
\underline{a}\neq 0]\\[0.2em]
\end{array}
\]}
\end{maplelatex}
\end{maplegroup}
\begin{maplegroup}
\begin{mapleinput}
\mapleinline{active}{1d}{Print(TD[2]);
}{}
\end{mapleinput}
\mapleresult
\begin{maplelatex}
\mapleinline{inert}{2d}{[a*xi1[x2]-xi3[x2] = 0, a^2*xi1[x3]+a*xi2[x2]-a[x2]*xi2+xi3[x1] = 0, a*xi1[x1]]-a*xi2[x2]+2*a[x2]*xi2-2*xi3[x1]-a*xi3[x3] = 0, a[x1] = 0, a[x2, x2] = 0, a[x2, x3] = 0, a[x3] = 0, a <> 0, xi2 <> 0]}{
\[
\begin{array}{l}
[a \, \underline{(\xi_1)_{x_2}}-(\xi_3)_{x_2}=0,\quad
a^2 \, \underline{(\xi_1)_{x_3}}+a \, (\xi_2)_{x_2}-a_{x_2} \, \xi_2+(\xi_3)_{x_1}=0,\\[0.5em] \qquad
a \, \underline{(\xi_1)_{x_1}}-a \, (\xi_2)_{x_2}+2\,a_{x_2} \, \xi_2-2\,(\xi_3)_{x_1}-a \, (\xi_3)_{x_3}=0,\quad
\underline{a_{x_1}}=0,\\[0.5em] \qquad
\underline{a_{x_2,x_2}}=0,\quad
\underline{a_{x_2,x_3}}=0,\quad
\underline{a_{x_3}}=0,\quad
\underline{a}\neq 0,\quad
\underline{\xi_2}\neq 0]\\[0.2em]
\end{array}
\]}
\end{maplelatex}
\end{maplegroup}
\begin{maplegroup}
\begin{mapleinput}
\mapleinline{active}{1d}{Print(TD[3]);
}{}
\end{mapleinput}
\mapleresult
\begin{maplelatex}
\mapleinline{inert}{2d}{[xi1[x1, x1]+xi2[x1, x2] = 0, xi1[x1, x2]+xi2[x2, x2] = 0, xi3[x1] = 0, xi3[x2] = 0, xi1[x1]+xi2[x2]-xi3[x3] = 0, a = 0]}{
\[
\begin{array}{l}
[\underline{(\xi_1)_{x_1,x_1}}+(\xi_2)_{x_1,x_2}=0,\quad
\underline{(\xi_1)_{x_1,x_2}}+(\xi_2)_{x_2,x_2}=0,\quad
\underline{(\xi_3)_{x_1}}=0,\quad
\underline{(\xi_3)_{x_2}}=0,\\[0.5em] \qquad
\underline{(\xi_1)_{x_1}}+(\xi_2)_{x_2}-(\xi_3)_{x_3}=0,\quad
\underline{a}=0]\\[0.2em]
\end{array}
\]}
\end{maplelatex}
\end{maplegroup}
\medskip
\noindent
With regard to the parametric function $a$, the first simple system is the
most generic one, in the sense that $a = a(x_2)$ is only assumed to be non-zero,
whereas in the second and third simple systems $a$ is subject to further equations.
In particular, the additional condition $a_{x_2,x_2} = 0$ derived
in \cite{PommaretQuadrat1997} to ensure formal integrability of the system
is exhibited in the second simple system of the Thomas decomposition.
\end{example}

\section{Conclusion}
\label{mlhdr:sec:conclusion}

In this paper the Thomas decomposition technique for systems of nonlinear
partial differential equations and inequations has been applied to
nonlinear control systems.
The method splits a given differential system into a finite family
of simple differential systems which are formally integrable and
define a partition of the solution set of the original differential system.
This symbolic approach allows to deal with both differential equations
and inequations, which may involve parameters.

Using elimination properties of the Thomas decomposition technique,
structural properties of nonlinear control systems have been investigated.
In particular, notions such as invertibility, observability and flat outputs
can be studied. In the presence of parameters, different simple systems
of a Thomas decomposition in general represent different structural behavior
of the control system.
A Maple implementation of Thomas' algorithm has been used to illustrate
the techniques on explicit examples.

At the time of this writing it is unclear how to adapt or generalize
the techniques to nonlinear differential time-delay systems or even
systems of nonlinear difference equations in full generality.
An analog of the notion of Thomas decomposition is not known for
systems of nonlinear difference equations.
However, note that, generalizing work of J.~F.\ Ritt \cite{Ritt} and
R.~M.\ Cohn \cite{Cohn}
and others, characteristic set methods have been developed for
ordinary difference polynomial systems and differential-difference
polynomial systems (cf., e.g., \cite{GaoLuoYuan}, \cite{GaovanderHoevenYuanZhang}).

\section*{acknowledgement}
The first author was partially supported by Schwerpunkt SPP~1489 of the
Deutsche Forschungsgemeinschaft.
The authors would like to thank an anonymous referee for several useful
remarks.
They would also like to thank S.~L.\ Rueda for pointing out reference~\cite{PicoMarco}.




\end{document}